\documentclass[11pt,a4paper]{article}
\usepackage[T1]{fontenc}
\usepackage{lmodern}
\input{glyphtounicode}
\usepackage{amsmath,amssymb,amsthm}
\usepackage{indentfirst}
\usepackage[a4paper,margin=25mm]{geometry}
\usepackage{graphicx}
\usepackage{tikz}
\usetikzlibrary{arrows.meta}
\definecolor{wtblue}{HTML}{2474CF}
\definecolor{wtblack}{HTML}{202020}
\definecolor{wtgray}{HTML}{7C8997}
\definecolor{wtpale}{HTML}{F5F8FD}
\tikzset{
  wt/base/.style={x=1cm,y=1cm,font=\small,line cap=round,line join=round},
  wt/edge/.style={draw=wtblack,line width=.70pt},
  wt/ghost/.style={draw=wtgray,line width=.70pt},
  wt/bracket/.style={draw=wtblue,line width=.55pt},
  wt/selection/.style={draw=wtblue,fill=wtpale,line width=.55pt,
    dash pattern=on 3pt off 2pt,rounded corners=2pt},
  wt/math/.style={inner sep=1pt,outer sep=0pt},
}

\usepackage{caption}
\usepackage{placeins}
\usepackage{microtype}
\usepackage[hidelinks,hypertexnames=false]{hyperref}
\hypersetup{pdftitle={Water Transport on 3-Spiders: Unbounded Lengths of Shortest Optimal Strategies},pdfauthor={Guanheng Chen, Ran Tao, Tianyi Tao}}
\title{\textbf{Water Transport on 3-Spiders:\\
Unbounded Lengths of Shortest Optimal Strategies}}

\newtheorem{theorem}{Theorem}
\newtheorem{lemma}{Lemma}

\theoremstyle{definition}
\newtheorem{definition}{Definition}

\date{}
\author{%
\begin{minipage}{.94\textwidth}
\centering
\normalsize
Guanheng Chen$^{1}$, Ran Tao$^{2}$, and Tianyi Tao$^{3}$\\[8pt]
\small
$^{1}$Shandong Normal University\\
\texttt{3577205810@qq.com}\\[5pt]
$^{2}$School of Mathematical Sciences and V.C.\ \& V.R.\ Key Lab,\\
Sichuan Normal University\\
\texttt{rtao\_sicnu@163.com}\\[5pt]
$^{3}$West China Second University Hospital, Sichuan University\\
\texttt{tytao@scu.edu.cn}
\end{minipage}}
\begin{document}
\maketitle
\begin{abstract}

Let $G$ be a finite simple graph. An averaging operation selects a nonempty
connected vertex set and replaces the weights of its vertices by their arithmetic
mean. For a nonnegative initial weight function $\omega:V(G)\to\mathbb R_{\ge0}$
and a target vertex $x$, let $t(G,x,\omega)$ denote the minimum number of
averaging operations needed to attain the maximum reachable weight at $x$.
Gollin et al.~\cite{gollin2024} proved, using the Axiom of Choice, that a finite
optimal strategy exists for every finite graph, every nonnegative initial weight
function and every target vertex; see Tao~\cite{tao2025} for a simplified
proof. We prove that the lengths of shortest optimal strategies on $3$-spiders
have no uniform finite bound. In particular, these lengths cannot be bounded
in terms of the maximum degree of the graph.

\end{abstract}
\section{The problem and its reformulation}

Let $G=(V,E)$ be a finite simple graph, let $\omega:V\to\mathbb R_{\ge0}$
be a nonnegative initial weight function, and fix a target vertex $x\in V$.
An \emph{averaging operation} selects a nonempty connected vertex set
$S\subseteq V$ and simultaneously replaces the weights of all vertices in $S$
by their arithmetic mean, leaving all other weights unchanged.
Let $t(G,x,\omega)$ be the minimum number of averaging operations needed to
attain the maximum reachable weight at $x$.

Water transport on graphs is closely related to averaging procedures for
consensus formation in the Deffuant model~\cite{haggstrom2012}.
Related problems on infinite and finite graphs were studied
in~\cite{haggstrom2019,vilkas2025}.
The existence of finite optimal strategies for connected-set averaging was
proved by Gollin et al.~\cite{gollin2024}; see also Tao~\cite{tao2025}
for a simplified proof and results on paths.

\begin{definition}
For an integer $m\ge2$, let $S(m,1,1)$ be the finite simple graph with vertex
set and edge set
\[
V\bigl(S(m,1,1)\bigr)
=\{a_0,a_1,\ldots,a_m,b,c\},
\]
\[
E\bigl(S(m,1,1)\bigr)
=
\bigl\{\{a_{i-1},a_i\}:1\le i\le m\bigr\}
\cup
\bigl\{\{a_0,b\},\{a_0,c\}\bigr\}
\]
where $a_0=o$. We call $S(m,1,1)$ a $3$-spider.
\end{definition}

More generally, write $\ell_i$ for the length of arm $i$ of a $3$-spider,
and $v_{i,j}$ for its $j$th vertex counting away from the root $o$.
We consider $\ell_1=m$ and $\ell_2=\ell_3=1$, where $m\ge2$, and write
$$
a_0=o,\qquad
a_j=v_{1,j}\ (1\le j\le m),\qquad
b=v_{2,1},\qquad c=v_{3,1}.
$$
The resulting graph is denoted by $S(m,1,1)$. Thus
$$
V\bigl(S(m,1,1)\bigr)
=
\{a_0,a_1,\ldots,a_m,b,c\},
$$
$$
E\bigl(S(m,1,1)\bigr)
=
\bigl\{\{a_{i-1},a_i\}:1\le i\le m\bigr\}
\cup
\bigl\{\{a_0,b\},\{a_0,c\}\bigr\}.
$$

To prove that the lengths of shortest optimal strategies on $3$-spiders have
no uniform finite bound, it suffices to show that, for every $m\ge2$, there
is a nonnegative weight function $\omega_m$ on $S(m,1,1)$ such that
$$
t\bigl(S(m,1,1),o,\omega_m\bigr)\ge m+2.
$$

\FloatBarrier
\section{Barrier lemmas}

Throughout this section, $S(m,1,1)$ has the vertex and edge sets specified
in Section~1. We use the same descriptions for $m=0,1$ when a shorter
remaining path occurs in Section~4. The proofs of Lemmas~1 and~2 apply
for $m\ge1$, and the proof of Lemma~3 applies for $m\ge0$.
Write $\omega(v)$ for the weight of a vertex $v$.
For $0\le k,t\le m+1$, define

$$
X_k=\sum_{i=0}^{k-1}\omega(a_i),\qquad
Y_t=\sum_{i=m-t+1}^{m}\omega(a_i),
$$

where

$$
X_0=Y_0=0.
$$

Here $X_k$ is the sum of the weights of the first $k$ vertices, starting
at $a_0=o$, of the path
$$
\langle a_0,a_1,\ldots,a_m\rangle
$$
and $Y_t$ is the sum of the weights of its last $t$ vertices.

\begin{definition}[Effective operation]
An averaging operation is \emph{effective} if the weights on the selected
connected vertex set are not all equal immediately before the operation.
\end{definition}

\begin{lemma}[Barrier for a source at the root]
Suppose that the initial weight function $\omega$ on $S(m,1,1)$ satisfies
$$
\omega(a_0)=1,\qquad
\omega(a_i)=0\quad(1\le i\le m),\qquad
\omega(b)=\omega(c)=0.
$$
Then, after any finite sequence of averaging operations,
$$
\omega(a_m)\le \frac{1}{m+1}.
$$

If equality holds, then the last effective operation containing $a_m$ must
average the set
$$
\{a_0,a_1,\ldots,a_m\}.
$$
Immediately after that operation,
$$
\omega(a_0)=\omega(a_1)=\cdots=\omega(a_m)=\frac{1}{m+1},
\qquad
\omega(b)=\omega(c)=0.
$$
\end{lemma}
\begin{proof}
For $1\le t\le m+1$, define
\[
T_t=\{a_{m-t+1},a_{m-t+2},\ldots,a_m\},
\]
and put $T_0=\varnothing$. Thus $T_t$ consists of the last $t$ vertices of
the path
\[
\langle a_0,a_1,\ldots,a_m\rangle
\]
counting from $a_m$ towards $a_0$. By the definition of $Y_t$,
\[
Y_t=\sum_{v\in T_t}\omega(v).
\]
The initial and terminal segments of the path, together with their vertex
counts, are shown in Figure~\ref{fig:path-sums}.

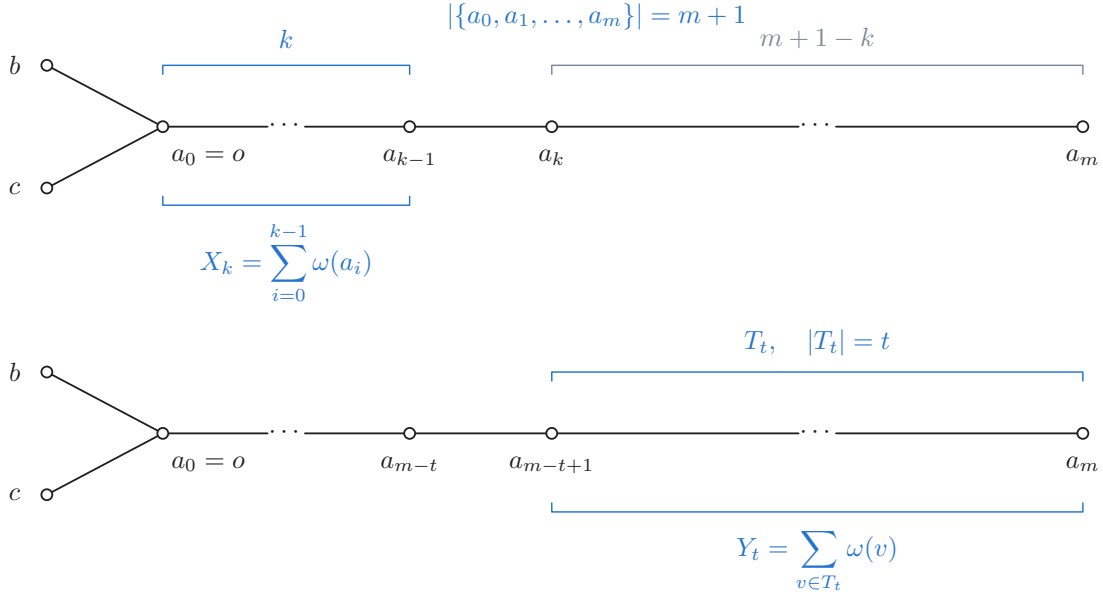
\begin{figure}[!htbp]
\centering
\resizebox{.96\linewidth}{!}{
\begin{tikzpicture}[wt/base]
\path (0,-6.4) rectangle (15.4,1.7);
\draw[wt/edge] (0.650,0.820) -- (2.200,0.000);
\draw[wt/edge] (0.650,-0.820) -- (2.200,0.000);
\draw[wt/edge] (2.200,0.000) -- (3.600,0.000);
\draw[wt/edge] (4.100,0.000) -- (5.500,0.000);
\draw[wt/edge] (5.500,0.000) -- (7.400,0.000);
\draw[wt/edge] (7.400,0.000) -- (10.700,0.000);
\draw[wt/edge] (11.200,0.000) -- (14.500,0.000);
\draw[wt/bracket,draw=wtblue] (2.200,0.700) -- (2.200,0.820) -- (5.500,0.820) -- (5.500,0.700);
\draw[wt/bracket,draw=wtgray] (7.400,0.700) -- (7.400,0.820) -- (14.500,0.820) -- (14.500,0.700);
\draw[wt/bracket,draw=wtblue] (2.200,-0.930) -- (2.200,-1.050) -- (5.500,-1.050) -- (5.500,-0.930);
\draw[wt/edge] (0.650,-3.280) -- (2.200,-4.100);
\draw[wt/edge] (0.650,-4.920) -- (2.200,-4.100);
\draw[wt/edge] (2.200,-4.100) -- (3.600,-4.100);
\draw[wt/edge] (4.100,-4.100) -- (5.500,-4.100);
\draw[wt/edge] (5.500,-4.100) -- (7.400,-4.100);
\draw[wt/edge] (7.400,-4.100) -- (10.700,-4.100);
\draw[wt/edge] (11.200,-4.100) -- (14.500,-4.100);
\draw[wt/bracket,draw=wtblue] (7.400,-3.400) -- (7.400,-3.280) -- (14.500,-3.280) -- (14.500,-3.400);
\draw[wt/bracket,draw=wtblue] (7.400,-5.030) -- (7.400,-5.150) -- (14.500,-5.150) -- (14.500,-5.030);
\draw[draw=wtblack,fill=white,line width=.70pt] (2.200,0.000) circle[radius=.073cm];
\draw[draw=wtblack,fill=white,line width=.70pt] (5.500,0.000) circle[radius=.073cm];
\draw[draw=wtblack,fill=white,line width=.70pt] (7.400,0.000) circle[radius=.073cm];
\draw[draw=wtblack,fill=white,line width=.70pt] (14.500,0.000) circle[radius=.073cm];
\draw[wt/edge,fill=white] (0.650,0.820) circle[radius=.073cm];
\draw[wt/edge,fill=white] (0.650,-0.820) circle[radius=.073cm];
\draw[draw=wtblack,fill=white,line width=.70pt] (2.200,-4.100) circle[radius=.073cm];
\draw[draw=wtblack,fill=white,line width=.70pt] (5.500,-4.100) circle[radius=.073cm];
\draw[draw=wtblack,fill=white,line width=.70pt] (7.400,-4.100) circle[radius=.073cm];
\draw[draw=wtblack,fill=white,line width=.70pt] (14.500,-4.100) circle[radius=.073cm];
\draw[wt/edge,fill=white] (0.650,-3.280) circle[radius=.073cm];
\draw[wt/edge,fill=white] (0.650,-4.920) circle[radius=.073cm];
\node[wt/math,text=wtblue,anchor=center] at (8.000,1.480) {$\displaystyle |\{a_0,a_1,\ldots,a_m\}|=m+1$};
\node[wt/math,text=wtblack,anchor=center] at (3.850,0.018) {$\displaystyle \cdots$};
\node[wt/math,text=wtblack,anchor=center] at (10.950,0.018) {$\displaystyle \cdots$};
\node[wt/math,text=wtblack,anchor=north west] at (2.280,-0.270) {$\displaystyle a_0=o$};
\node[wt/math,text=wtblack,anchor=north] at (5.500,-0.270) {$\displaystyle a_{k-1}$};
\node[wt/math,text=wtblack,anchor=north] at (7.400,-0.270) {$\displaystyle a_k$};
\node[wt/math,text=wtblack,anchor=north] at (14.500,-0.270) {$\displaystyle a_m$};
\node[wt/math,text=wtblack,anchor=east] at (0.330,0.820) {$\displaystyle b$};
\node[wt/math,text=wtblack,anchor=east] at (0.330,-0.820) {$\displaystyle c$};
\node[wt/math,text=wtblue,anchor=south] at (3.850,1.010) {$\displaystyle k$};
\node[wt/math,text=wtgray,anchor=south] at (10.950,1.010) {$\displaystyle m+1-k$};
\node[wt/math,text=wtblue,anchor=north] at (3.850,-1.240) {$\displaystyle X_k=\sum_{i=0}^{k-1}\omega(a_i)$};
\node[wt/math,text=wtblack,anchor=center] at (3.850,-4.082) {$\displaystyle \cdots$};
\node[wt/math,text=wtblack,anchor=center] at (10.950,-4.082) {$\displaystyle \cdots$};
\node[wt/math,text=wtblack,anchor=north west] at (2.280,-4.370) {$\displaystyle a_0=o$};
\node[wt/math,text=wtblack,anchor=north] at (5.500,-4.370) {$\displaystyle a_{m-t}$};
\node[wt/math,text=wtblack,anchor=north] at (7.400,-4.370) {$\displaystyle a_{m-t+1}$};
\node[wt/math,text=wtblack,anchor=north] at (14.500,-4.370) {$\displaystyle a_m$};
\node[wt/math,text=wtblack,anchor=east] at (0.330,-3.280) {$\displaystyle b$};
\node[wt/math,text=wtblack,anchor=east] at (0.330,-4.920) {$\displaystyle c$};
\node[wt/math,text=wtblue,anchor=south] at (10.950,-3.090) {$\displaystyle T_t,\quad |T_t|=t$};
\node[wt/math,text=wtblue,anchor=north] at (10.950,-5.340) {$\displaystyle Y_t=\sum_{v\in T_t}\omega(v)$};
\end{tikzpicture}
}
\caption{Prefix sums and terminal sums along the path.}
\label{fig:path-sums}
\end{figure}

Averaging preserves nonnegativity and the sum of all vertex weights.
Hence, throughout the sequence,
\[
\omega(v)\geq 0
\quad
\bigl(v\in V(S(m,1,1))\bigr),
\qquad
\sum_{v\in V(S(m,1,1))}\omega(v)=1.
\]

We prove by induction on the number of operations that
\begin{equation}\label{eq:root-tail-bound}
Y_t\leq \frac{t}{m+1},
\qquad 0\leq t\leq m+1,
\end{equation}
and, simultaneously, that for $1\le t\le m$,
\begin{equation}\label{eq:root-tail-rigidity}
Y_t=\frac{t}{m+1}
\quad\Longrightarrow\quad
\omega(a_{m-t+1})=\cdots=\omega(a_m)=\frac{1}{m+1}.
\end{equation}
Taking $t=1$ in \eqref{eq:root-tail-bound} then gives
\[
\omega(a_m)=Y_1\leq \frac{1}{m+1},
\]
which is the required bound.

Initially, $Y_t=0$ for $0\le t\le m$, and
\[
Y_{m+1}=\sum_{i=0}^{m}\omega(a_i)=1.
\]
Thus \eqref{eq:root-tail-bound} holds. For $1\le t\le m$, the equality $Y_t=t/(m+1)$ does not
hold initially, so \eqref{eq:root-tail-rigidity} also holds.

Suppose that the current weight function $\omega$ satisfies \eqref{eq:root-tail-bound} and \eqref{eq:root-tail-rigidity}.
Let $S$ be the connected set selected by the next operation, and let $\omega'$
be the resulting weight function. For $1\le t\le m+1$, put
\[
Y_t'=\sum_{v\in T_t}\omega'(v)
=\sum_{i=m-t+1}^{m}\omega'(a_i),
\]
and put $Y_0'=0$. Thus $Y_t'$ is the total weight on $T_t$ after the
operation. We show that $\omega'$ also satisfies \eqref{eq:root-tail-bound} and \eqref{eq:root-tail-rigidity}.

If the operation is ineffective, then $\omega'=\omega$ and there is nothing to prove.
We may therefore assume that it is effective.

First suppose that
\[
S\cap\{b,c\}=\varnothing.
\]
Since $S$ is connected, there are integers $0\le r\le j\le m$ such that
\[
S=\{a_r,a_{r+1},\ldots,a_j\}.
\]
Define
\[
t_-=m-j,
\qquad
t_+=m-r+1.
\]
Thus $t_-$ and $t_+$ count the vertices of $\{a_{j+1},\ldots,a_m\}$ and
$\{a_r,\ldots,a_m\}$, respectively. These numbers distinguish the three
possible positions of $T_t$ relative to $S$:
\[
t\leq t_-
\iff T_t\cap S=\varnothing,
\]
\[
t\geq t_+
\iff S\subseteq T_t,
\]
and
\[
t_-<t<t_+
\iff T_t\cap S\neq\varnothing
\ \text{ and }\ 
S\not\subseteq T_t.
\]
See Figure~\ref{fig:three-positions}.

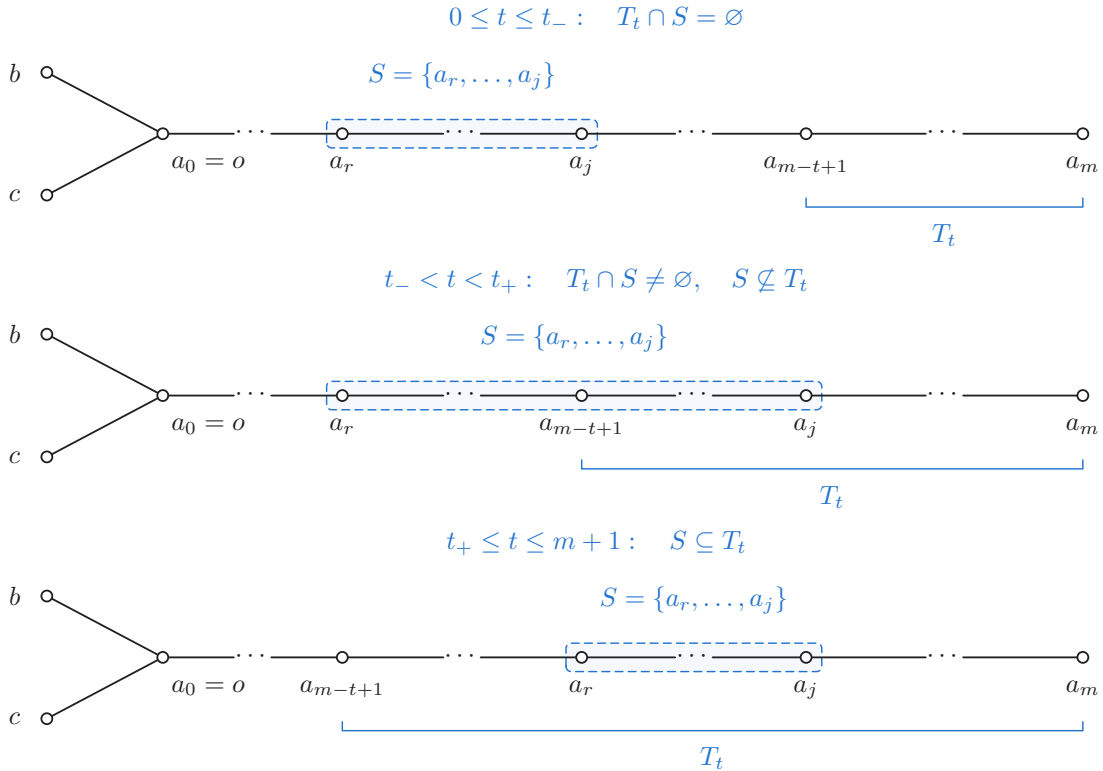
\begin{figure}[!htbp]
\centering
\resizebox{.96\linewidth}{!}{
\begin{tikzpicture}[wt/base]
\path (0,-8.7) rectangle (15.4,1.8);
\draw[wt/selection] (4.390,-0.190) rectangle (8.010,0.190);
\draw[wt/selection] (4.390,-3.690) rectangle (11.010,-3.310);
\draw[wt/selection] (7.590,-7.190) rectangle (11.010,-6.810);
\draw[wt/edge] (0.650,0.820) -- (2.200,0.000);
\draw[wt/edge] (0.650,-0.820) -- (2.200,0.000);
\draw[wt/edge] (2.200,0.000) -- (3.150,0.000);
\draw[wt/edge] (3.650,0.000) -- (4.600,0.000);
\draw[wt/edge] (4.600,0.000) -- (5.950,0.000);
\draw[wt/edge] (6.450,0.000) -- (7.800,0.000);
\draw[wt/edge] (7.800,0.000) -- (9.050,0.000);
\draw[wt/edge] (9.550,0.000) -- (10.800,0.000);
\draw[wt/edge] (10.800,0.000) -- (12.400,0.000);
\draw[wt/edge] (12.900,0.000) -- (14.500,0.000);
\draw[wt/bracket,draw=wtblue] (10.800,-0.860) -- (10.800,-0.980) -- (14.500,-0.980) -- (14.500,-0.860);
\draw[wt/edge] (0.650,-2.680) -- (2.200,-3.500);
\draw[wt/edge] (0.650,-4.320) -- (2.200,-3.500);
\draw[wt/edge] (2.200,-3.500) -- (3.150,-3.500);
\draw[wt/edge] (3.650,-3.500) -- (4.600,-3.500);
\draw[wt/edge] (4.600,-3.500) -- (5.950,-3.500);
\draw[wt/edge] (6.450,-3.500) -- (7.800,-3.500);
\draw[wt/edge] (7.800,-3.500) -- (9.050,-3.500);
\draw[wt/edge] (9.550,-3.500) -- (10.800,-3.500);
\draw[wt/edge] (10.800,-3.500) -- (12.400,-3.500);
\draw[wt/edge] (12.900,-3.500) -- (14.500,-3.500);
\draw[wt/bracket,draw=wtblue] (7.800,-4.360) -- (7.800,-4.480) -- (14.500,-4.480) -- (14.500,-4.360);
\draw[wt/edge] (0.650,-6.180) -- (2.200,-7.000);
\draw[wt/edge] (0.650,-7.820) -- (2.200,-7.000);
\draw[wt/edge] (2.200,-7.000) -- (3.150,-7.000);
\draw[wt/edge] (3.650,-7.000) -- (4.600,-7.000);
\draw[wt/edge] (4.600,-7.000) -- (5.950,-7.000);
\draw[wt/edge] (6.450,-7.000) -- (7.800,-7.000);
\draw[wt/edge] (7.800,-7.000) -- (9.050,-7.000);
\draw[wt/edge] (9.550,-7.000) -- (10.800,-7.000);
\draw[wt/edge] (10.800,-7.000) -- (12.400,-7.000);
\draw[wt/edge] (12.900,-7.000) -- (14.500,-7.000);
\draw[wt/bracket,draw=wtblue] (4.600,-7.860) -- (4.600,-7.980) -- (14.500,-7.980) -- (14.500,-7.860);
\draw[draw=wtblack,fill=white,line width=.70pt] (2.200,0.000) circle[radius=.073cm];
\draw[draw=wtblack,fill=white,line width=.70pt] (4.600,0.000) circle[radius=.073cm];
\draw[draw=wtblack,fill=white,line width=.70pt] (7.800,0.000) circle[radius=.073cm];
\draw[draw=wtblack,fill=white,line width=.70pt] (10.800,0.000) circle[radius=.073cm];
\draw[draw=wtblack,fill=white,line width=.70pt] (14.500,0.000) circle[radius=.073cm];
\draw[wt/edge,fill=white] (0.650,0.820) circle[radius=.073cm];
\draw[wt/edge,fill=white] (0.650,-0.820) circle[radius=.073cm];
\draw[draw=wtblack,fill=white,line width=.70pt] (2.200,-3.500) circle[radius=.073cm];
\draw[draw=wtblack,fill=white,line width=.70pt] (4.600,-3.500) circle[radius=.073cm];
\draw[draw=wtblack,fill=white,line width=.70pt] (7.800,-3.500) circle[radius=.073cm];
\draw[draw=wtblack,fill=white,line width=.70pt] (10.800,-3.500) circle[radius=.073cm];
\draw[draw=wtblack,fill=white,line width=.70pt] (14.500,-3.500) circle[radius=.073cm];
\draw[wt/edge,fill=white] (0.650,-2.680) circle[radius=.073cm];
\draw[wt/edge,fill=white] (0.650,-4.320) circle[radius=.073cm];
\draw[draw=wtblack,fill=white,line width=.70pt] (2.200,-7.000) circle[radius=.073cm];
\draw[draw=wtblack,fill=white,line width=.70pt] (4.600,-7.000) circle[radius=.073cm];
\draw[draw=wtblack,fill=white,line width=.70pt] (7.800,-7.000) circle[radius=.073cm];
\draw[draw=wtblack,fill=white,line width=.70pt] (10.800,-7.000) circle[radius=.073cm];
\draw[draw=wtblack,fill=white,line width=.70pt] (14.500,-7.000) circle[radius=.073cm];
\draw[wt/edge,fill=white] (0.650,-6.180) circle[radius=.073cm];
\draw[wt/edge,fill=white] (0.650,-7.820) circle[radius=.073cm];
\node[wt/math,text=wtblack,anchor=center] at (3.400,0.018) {$\displaystyle \cdots$};
\node[wt/math,text=wtblack,anchor=center] at (6.200,0.018) {$\displaystyle \cdots$};
\node[wt/math,text=wtblack,anchor=center] at (9.300,0.018) {$\displaystyle \cdots$};
\node[wt/math,text=wtblack,anchor=center] at (12.650,0.018) {$\displaystyle \cdots$};
\node[wt/math,text=wtblack,anchor=north west] at (2.280,-0.270) {$\displaystyle a_0=o$};
\node[wt/math,text=wtblack,anchor=north] at (4.600,-0.270) {$\displaystyle a_r$};
\node[wt/math,text=wtblack,anchor=north] at (7.800,-0.270) {$\displaystyle a_j$};
\node[wt/math,text=wtblack,anchor=north] at (10.800,-0.270) {$\displaystyle a_{m-t+1}$};
\node[wt/math,text=wtblack,anchor=north] at (14.500,-0.270) {$\displaystyle a_m$};
\node[wt/math,text=wtblack,anchor=east] at (0.330,0.820) {$\displaystyle b$};
\node[wt/math,text=wtblack,anchor=east] at (0.330,-0.820) {$\displaystyle c$};
\node[wt/math,text=wtblue,anchor=center] at (8.000,1.530) {$\displaystyle 0\le t\le t_-:\quad T_t\cap S=\varnothing$};
\node[wt/math,text=wtblue,anchor=south] at (6.200,0.550) {$\displaystyle S=\{a_r,\ldots,a_j\}$};
\node[wt/math,text=wtblue,anchor=north] at (12.650,-1.170) {$\displaystyle T_t$};
\node[wt/math,text=wtblack,anchor=center] at (3.400,-3.482) {$\displaystyle \cdots$};
\node[wt/math,text=wtblack,anchor=center] at (6.200,-3.482) {$\displaystyle \cdots$};
\node[wt/math,text=wtblack,anchor=center] at (9.300,-3.482) {$\displaystyle \cdots$};
\node[wt/math,text=wtblack,anchor=center] at (12.650,-3.482) {$\displaystyle \cdots$};
\node[wt/math,text=wtblack,anchor=north west] at (2.280,-3.770) {$\displaystyle a_0=o$};
\node[wt/math,text=wtblack,anchor=north] at (4.600,-3.770) {$\displaystyle a_r$};
\node[wt/math,text=wtblack,anchor=north] at (7.800,-3.770) {$\displaystyle a_{m-t+1}$};
\node[wt/math,text=wtblack,anchor=north] at (10.800,-3.770) {$\displaystyle a_j$};
\node[wt/math,text=wtblack,anchor=north] at (14.500,-3.770) {$\displaystyle a_m$};
\node[wt/math,text=wtblack,anchor=east] at (0.330,-2.680) {$\displaystyle b$};
\node[wt/math,text=wtblack,anchor=east] at (0.330,-4.320) {$\displaystyle c$};
\node[wt/math,text=wtblue,anchor=center] at (8.000,-1.970) {$\displaystyle t_-<t<t_+:\quad T_t\cap S\ne\varnothing,\quad S\not\subseteq T_t$};
\node[wt/math,text=wtblue,anchor=south] at (7.700,-2.950) {$\displaystyle S=\{a_r,\ldots,a_j\}$};
\node[wt/math,text=wtblue,anchor=north] at (11.150,-4.670) {$\displaystyle T_t$};
\node[wt/math,text=wtblack,anchor=center] at (3.400,-6.982) {$\displaystyle \cdots$};
\node[wt/math,text=wtblack,anchor=center] at (6.200,-6.982) {$\displaystyle \cdots$};
\node[wt/math,text=wtblack,anchor=center] at (9.300,-6.982) {$\displaystyle \cdots$};
\node[wt/math,text=wtblack,anchor=center] at (12.650,-6.982) {$\displaystyle \cdots$};
\node[wt/math,text=wtblack,anchor=north west] at (2.280,-7.270) {$\displaystyle a_0=o$};
\node[wt/math,text=wtblack,anchor=north] at (4.600,-7.270) {$\displaystyle a_{m-t+1}$};
\node[wt/math,text=wtblack,anchor=north] at (7.800,-7.270) {$\displaystyle a_r$};
\node[wt/math,text=wtblack,anchor=north] at (10.800,-7.270) {$\displaystyle a_j$};
\node[wt/math,text=wtblack,anchor=north] at (14.500,-7.270) {$\displaystyle a_m$};
\node[wt/math,text=wtblack,anchor=east] at (0.330,-6.180) {$\displaystyle b$};
\node[wt/math,text=wtblack,anchor=east] at (0.330,-7.820) {$\displaystyle c$};
\node[wt/math,text=wtblue,anchor=center] at (8.000,-5.470) {$\displaystyle t_+\le t\le m+1:\quad S\subseteq T_t$};
\node[wt/math,text=wtblue,anchor=south] at (9.300,-6.450) {$\displaystyle S=\{a_r,\ldots,a_j\}$};
\node[wt/math,text=wtblue,anchor=north] at (9.550,-8.170) {$\displaystyle T_t$};
\end{tikzpicture}
}
\caption{The three possible positions of $T_t$ relative to the path interval $S$.}
\label{fig:three-positions}
\end{figure}

We consider these three cases in turn.

If $0\le t\le t_-$, then $T_t\cap S=\varnothing$. No weight on $T_t$
changes, and hence
\[
Y_t'=Y_t\leq \frac{t}{m+1}.
\]

If $t_+\le t\le m+1$, then $S\subseteq T_t$. The operation preserves the
sum of the weights on $S$ and leaves all weights on $T_t\setminus S$
unchanged. Thus
\[
Y_t'=Y_t\leq \frac{t}{m+1}.
\]

Now suppose that
\[
t_-<t<t_+.
\]
Then
\[
T_t\setminus S
=\{a_{j+1},a_{j+2},\ldots,a_m\}
=T_{t_-},
\]
and
\[
T_t\cap S
=\{a_{m-t+1},a_{m-t+2},\ldots,a_j\}.
\]
The intersection $T_t\cap S$ contains $t-t_-$ vertices.

Before the operation, the sum of the weights on $S$ is
\[
\sum_{i=r}^{j}\omega(a_i)=Y_{t_+}-Y_{t_-},
\]
and $S$ contains $t_+-t_-$ vertices. If $\mu$ denotes their common weight
after the operation, then
\[
\mu=\frac{Y_{t_+}-Y_{t_-}}{t_+-t_-}.
\]
The weights on $T_{t_-}$ do not change, while all $t-t_-$ vertices of
$T_t\cap S$ receive weight $\mu$. Therefore,
\begin{equation}\label{eq:path-interpolation}
\begin{aligned}
Y_t'
&=Y_{t_-}+(t-t_-)\mu\\
&=\frac{t_+-t}{t_+-t_-}Y_{t_-}
  +\frac{t-t_-}{t_+-t_-}Y_{t_+}\\
&\leq
\frac{t_+-t}{t_+-t_-}\frac{t_-}{m+1}
+\frac{t-t_-}{t_+-t_-}\frac{t_+}{m+1}\\
&=\frac{t}{m+1}.
\end{aligned}
\end{equation}
The vertex counts in \eqref{eq:path-interpolation} are illustrated in
Figure~\ref{fig:interpolation-counts}.

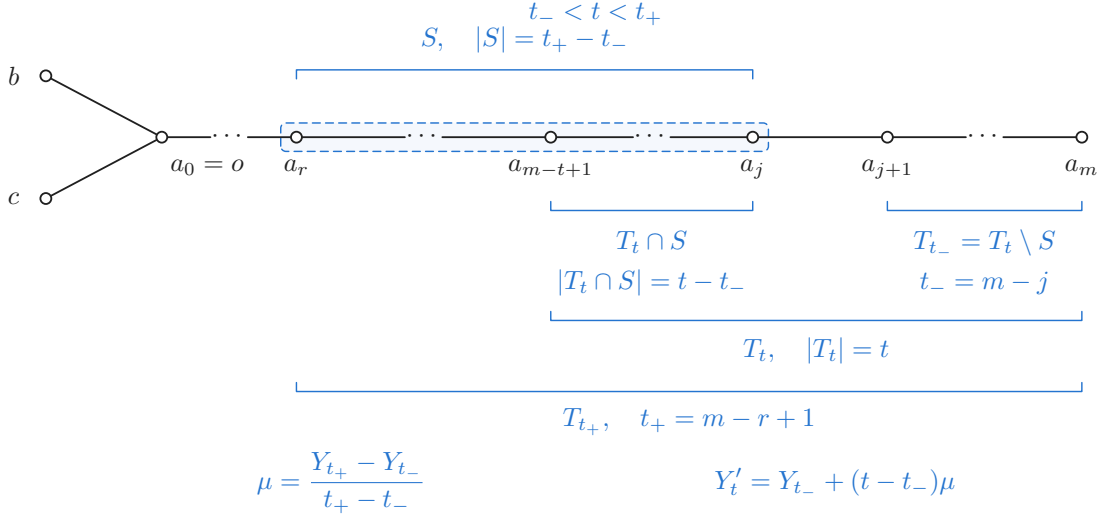
\begin{figure}[!htbp]
\centering
\resizebox{.96\linewidth}{!}{
\begin{tikzpicture}[wt/base]
\path (0,-5.3) rectangle (15.4,1.9);
\draw[wt/selection] (3.790,-0.190) rectangle (10.310,0.190);
\draw[wt/edge] (0.650,0.820) -- (2.200,0.000);
\draw[wt/edge] (0.650,-0.820) -- (2.200,0.000);
\draw[wt/edge] (2.200,0.000) -- (2.850,0.000);
\draw[wt/edge] (3.350,0.000) -- (4.000,0.000);
\draw[wt/edge] (4.000,0.000) -- (5.450,0.000);
\draw[wt/edge] (5.950,0.000) -- (7.400,0.000);
\draw[wt/edge] (7.400,0.000) -- (8.500,0.000);
\draw[wt/edge] (9.000,0.000) -- (10.100,0.000);
\draw[wt/edge] (10.100,0.000) -- (11.900,0.000);
\draw[wt/edge] (11.900,0.000) -- (12.950,0.000);
\draw[wt/edge] (13.450,0.000) -- (14.500,0.000);
\draw[wt/bracket,draw=wtblue] (4.000,0.780) -- (4.000,0.900) -- (10.100,0.900) -- (10.100,0.780);
\draw[wt/bracket,draw=wtblue] (7.400,-0.860) -- (7.400,-0.980) -- (10.100,-0.980) -- (10.100,-0.860);
\draw[wt/bracket,draw=wtblue] (11.900,-0.860) -- (11.900,-0.980) -- (14.500,-0.980) -- (14.500,-0.860);
\draw[wt/bracket,draw=wtblue] (7.400,-2.330) -- (7.400,-2.450) -- (14.500,-2.450) -- (14.500,-2.330);
\draw[wt/bracket,draw=wtblue] (4.000,-3.280) -- (4.000,-3.400) -- (14.500,-3.400) -- (14.500,-3.280);
\draw[draw=wtblack,fill=white,line width=.70pt] (2.200,0.000) circle[radius=.073cm];
\draw[draw=wtblack,fill=white,line width=.70pt] (4.000,0.000) circle[radius=.073cm];
\draw[draw=wtblack,fill=white,line width=.70pt] (7.400,0.000) circle[radius=.073cm];
\draw[draw=wtblack,fill=white,line width=.70pt] (10.100,0.000) circle[radius=.073cm];
\draw[draw=wtblack,fill=white,line width=.70pt] (11.900,0.000) circle[radius=.073cm];
\draw[draw=wtblack,fill=white,line width=.70pt] (14.500,0.000) circle[radius=.073cm];
\draw[wt/edge,fill=white] (0.650,0.820) circle[radius=.073cm];
\draw[wt/edge,fill=white] (0.650,-0.820) circle[radius=.073cm];
\node[wt/math,text=wtblack,anchor=center] at (3.100,0.018) {$\displaystyle \cdots$};
\node[wt/math,text=wtblack,anchor=center] at (5.700,0.018) {$\displaystyle \cdots$};
\node[wt/math,text=wtblack,anchor=center] at (8.750,0.018) {$\displaystyle \cdots$};
\node[wt/math,text=wtblack,anchor=center] at (13.200,0.018) {$\displaystyle \cdots$};
\node[wt/math,text=wtblack,anchor=north west] at (2.280,-0.270) {$\displaystyle a_0=o$};
\node[wt/math,text=wtblack,anchor=north] at (4.000,-0.270) {$\displaystyle a_r$};
\node[wt/math,text=wtblack,anchor=north] at (7.400,-0.270) {$\displaystyle a_{m-t+1}$};
\node[wt/math,text=wtblack,anchor=north] at (10.100,-0.270) {$\displaystyle a_j$};
\node[wt/math,text=wtblack,anchor=north] at (11.900,-0.270) {$\displaystyle a_{j+1}$};
\node[wt/math,text=wtblack,anchor=north] at (14.500,-0.270) {$\displaystyle a_m$};
\node[wt/math,text=wtblack,anchor=east] at (0.330,0.820) {$\displaystyle b$};
\node[wt/math,text=wtblack,anchor=east] at (0.330,-0.820) {$\displaystyle c$};
\node[wt/math,text=wtblue,anchor=center] at (8.000,1.650) {$\displaystyle t_-<t<t_+$};
\node[wt/math,text=wtblue,anchor=south] at (7.050,1.090) {$\displaystyle S,\quad |S|=t_+-t_-$};
\node[wt/math,text=wtblue,anchor=north] at (8.750,-1.170) {$\displaystyle \begin{gathered}T_t\cap S\\ |T_t\cap S|=t-t_-\end{gathered}$};
\node[wt/math,text=wtblue,anchor=north] at (13.200,-1.170) {$\displaystyle \begin{gathered}T_{t_-}=T_t\setminus S\\ t_-=m-j\end{gathered}$};
\node[wt/math,text=wtblue,anchor=north] at (10.950,-2.640) {$\displaystyle T_t,\quad |T_t|=t$};
\node[wt/math,text=wtblue,anchor=north] at (9.250,-3.590) {$\displaystyle T_{t_+},\quad t_+=m-r+1$};
\node[wt/math,text=wtblue,anchor=center] at (4.600,-4.620) {$\displaystyle \mu=\frac{Y_{t_+}-Y_{t_-}}{t_+-t_-}$};
\node[wt/math,text=wtblue,anchor=center] at (11.200,-4.620) {$\displaystyle Y_t'=Y_{t_-}+(t-t_-)\mu$};
\end{tikzpicture}
}
\caption{Set decomposition and vertex counts in the interpolation identity.}
\label{fig:interpolation-counts}
\end{figure}

This proves that $\omega'$ satisfies \eqref{eq:root-tail-bound}.

We next prove \eqref{eq:root-tail-rigidity} for $\omega'$. Suppose that $1\le t\le m$ and
\[
Y_t'=\frac{t}{m+1}.
\]

If $t\le t_-$, then $Y_t'=Y_t$ and $T_t\cap S=\varnothing$. By \eqref{eq:root-tail-rigidity},
before the operation,
\[
\omega(a_{m-t+1})=\cdots=\omega(a_m)=\frac{1}{m+1}.
\]
Since no weight on $T_t$ changes, after the operation we have
\[
\omega'(a_{m-t+1})=\cdots=\omega'(a_m)=\frac{1}{m+1}.
\]

If $t\ge t_+$, then $Y_t'=Y_t$ and $S\subseteq T_t$. By \eqref{eq:root-tail-rigidity}, every
vertex of $T_t$, and hence every vertex of $S$, has weight $1/(m+1)$
before the operation. Averaging $S$ leaves these weights unchanged.
The vertices of $T_t\setminus S$ are not selected, so
\[
\omega'(a_{m-t+1})=\cdots=\omega'(a_m)=\frac{1}{m+1}.
\]

Finally, suppose that $t_-<t<t_+$. The coefficients of $Y_{t_-}$ and
$Y_{t_+}$ in \eqref{eq:path-interpolation} are
\[
\frac{t_+-t}{t_+-t_-}
\qquad\text{ and }\qquad
\frac{t-t_-}{t_+-t_-}.
\]
Both coefficients are positive, and
\[
\frac{t_+-t}{t_+-t_-}
+
\frac{t-t_-}{t_+-t_-}
=1.
\]
Consequently, equality in the final bound in \eqref{eq:path-interpolation} forces equality in \eqref{eq:root-tail-bound}
at both $t_-$ and $t_+$:
\[
Y_{t_-}=\frac{t_-}{m+1},
\qquad
Y_{t_+}=\frac{t_+}{m+1}.
\]
It follows that
\[
\mu
=\frac{Y_{t_+}-Y_{t_-}}{t_+-t_-}
=\frac{1}{m+1}.
\]

If $t_->0$, then \eqref{eq:root-tail-rigidity} gives, before the operation,
\[
\omega(a_{j+1})=\cdots=\omega(a_m)=\frac{1}{m+1}.
\]
These vertices form the set
\[
T_{t_-}=\{a_{j+1},\ldots,a_m\},
\]
and none is selected, so their weights do not change. If $t_-=0$, then
$T_{t_-}=\varnothing$. On the other hand,
\[
T_t\cap S
=\{a_{m-t+1},\ldots,a_j\},
\]
and every vertex of this intersection has weight
\[
\mu=\frac{1}{m+1}
\]
after the operation. Since
\[
T_t=T_{t_-}\cup(T_t\cap S),
\]
we obtain
\[
\omega'(a_{m-t+1})=\cdots=\omega'(a_m)=\frac{1}{m+1}.
\]
Thus $\omega'$ satisfies both \eqref{eq:root-tail-bound} and \eqref{eq:root-tail-rigidity} when $S\cap\{b,c\}=\varnothing$.

Next suppose that
\[
S\cap\{b,c\}\neq\varnothing.
\]
Define
\[
F=S\cap\{b,c\},
\qquad
s=|F|.
\]
Thus $F$ is the set of selected vertices among $b,c$, and $s$ is their
number. Since $F\ne\varnothing$,
\[
s\in\{1,2\}.
\]

Define
\[
k=\left|S\cap\{a_0,a_1,\ldots,a_m\}\right|,
\]
so that $k$ is the number of selected vertices on the path
\[
\langle a_0,a_1,\ldots,a_m\rangle.
\]
Since the operation is effective, $S$ cannot consist of the single vertex
$b$ or $c$. Connectivity therefore gives
\[
1\leq k\leq m+1
\]
and
\[
S\cap\{a_0,a_1,\ldots,a_m\}
=\{a_0,a_1,\ldots,a_{k-1}\}.
\]
Hence
\[
S=\{a_0,a_1,\ldots,a_{k-1}\}\cup F,
\qquad
|S|=k+s.
\]
The possible shapes of $S$ are shown in Figure~\ref{fig:selected-leaves}.

\begin{figure}[!htbp]
\centering
\resizebox{.96\linewidth}{!}{
\begin{tikzpicture}[wt/base]
\path (0,-6.8) rectangle (15.4,1.8);
\draw[wt/selection] (0.460,1.020) -- (0.840,1.020) -- (2.300,0.190) -- (6.910,0.190) -- (6.910,-0.190) -- (2.100,-0.190) -- (0.460,0.630) -- cycle;
\draw[wt/selection] (0.460,-3.280) -- (0.840,-3.280) -- (2.300,-4.110) -- (6.910,-4.110) -- (6.910,-4.490) -- (2.300,-4.490) -- (0.840,-5.320) -- (0.460,-5.320) -- cycle;
\draw[wt/edge] (0.650,0.820) -- (2.200,0.000);
\draw[wt/edge] (0.650,-0.820) -- (2.200,0.000);
\draw[wt/edge] (2.200,0.000) -- (4.200,0.000);
\draw[wt/edge] (4.700,0.000) -- (6.700,0.000);
\draw[wt/edge] (6.700,0.000) -- (8.400,0.000);
\draw[wt/edge] (8.400,0.000) -- (11.200,0.000);
\draw[wt/edge] (11.700,0.000) -- (14.500,0.000);
\draw[wt/bracket,draw=wtblue] (2.200,-0.900) -- (2.200,-1.020) -- (6.700,-1.020) -- (6.700,-0.900);
\draw[wt/bracket,draw=wtblue] (8.400,-0.900) -- (8.400,-1.020) -- (14.500,-1.020) -- (14.500,-0.900);
\draw[wt/edge] (0.650,-3.480) -- (2.200,-4.300);
\draw[wt/edge] (0.650,-5.120) -- (2.200,-4.300);
\draw[wt/edge] (2.200,-4.300) -- (4.200,-4.300);
\draw[wt/edge] (4.700,-4.300) -- (6.700,-4.300);
\draw[wt/edge] (6.700,-4.300) -- (8.400,-4.300);
\draw[wt/edge] (8.400,-4.300) -- (11.200,-4.300);
\draw[wt/edge] (11.700,-4.300) -- (14.500,-4.300);
\draw[wt/bracket,draw=wtblue] (2.200,-5.200) -- (2.200,-5.320) -- (6.700,-5.320) -- (6.700,-5.200);
\draw[wt/bracket,draw=wtblue] (8.400,-5.200) -- (8.400,-5.320) -- (14.500,-5.320) -- (14.500,-5.200);
\draw[draw=wtblack,fill=white,line width=.70pt] (2.200,0.000) circle[radius=.073cm];
\draw[draw=wtblack,fill=white,line width=.70pt] (6.700,0.000) circle[radius=.073cm];
\draw[draw=wtblack,fill=white,line width=.70pt] (8.400,0.000) circle[radius=.073cm];
\draw[draw=wtblack,fill=white,line width=.70pt] (14.500,0.000) circle[radius=.073cm];
\draw[wt/edge,fill=white] (0.650,0.820) circle[radius=.073cm];
\draw[wt/edge,fill=white] (0.650,-0.820) circle[radius=.073cm];
\draw[draw=wtblack,fill=white,line width=.70pt] (2.200,-4.300) circle[radius=.073cm];
\draw[draw=wtblack,fill=white,line width=.70pt] (6.700,-4.300) circle[radius=.073cm];
\draw[draw=wtblack,fill=white,line width=.70pt] (8.400,-4.300) circle[radius=.073cm];
\draw[draw=wtblack,fill=white,line width=.70pt] (14.500,-4.300) circle[radius=.073cm];
\draw[wt/edge,fill=white] (0.650,-3.480) circle[radius=.073cm];
\draw[wt/edge,fill=white] (0.650,-5.120) circle[radius=.073cm];
\node[wt/math,text=wtblack,anchor=center] at (4.450,0.018) {$\displaystyle \cdots$};
\node[wt/math,text=wtblack,anchor=center] at (11.450,0.018) {$\displaystyle \cdots$};
\node[wt/math,text=wtblack,anchor=north west] at (2.280,-0.270) {$\displaystyle a_0=o$};
\node[wt/math,text=wtblack,anchor=north] at (6.700,-0.270) {$\displaystyle a_{k-1}$};
\node[wt/math,text=wtblack,anchor=north] at (8.400,-0.270) {$\displaystyle a_k$};
\node[wt/math,text=wtblack,anchor=north] at (14.500,-0.270) {$\displaystyle a_m$};
\node[wt/math,text=wtblack,anchor=east] at (0.330,0.820) {$\displaystyle b$};
\node[wt/math,text=wtblack,anchor=east] at (0.330,-0.820) {$\displaystyle c$};
\node[wt/math,text=wtblue,anchor=center] at (8.000,1.530) {$\displaystyle F=\{b\},\quad s=1$};
\node[wt/math,text=wtblue,anchor=south] at (4.600,0.470) {$\displaystyle S$};
\node[wt/math,text=wtblue,anchor=north] at (4.450,-1.210) {$\displaystyle k$};
\node[wt/math,text=wtblue,anchor=north] at (11.450,-1.210) {$\displaystyle m+1-k$};
\node[wt/math,text=wtblue,anchor=center] at (8.000,-1.920) {$\displaystyle S=\{a_0,\ldots,a_{k-1}\}\cup F,\quad |S|=k+s$};
\node[wt/math,text=wtblack,anchor=center] at (4.450,-4.282) {$\displaystyle \cdots$};
\node[wt/math,text=wtblack,anchor=center] at (11.450,-4.282) {$\displaystyle \cdots$};
\node[wt/math,text=wtblack,anchor=north west] at (2.280,-4.570) {$\displaystyle a_0=o$};
\node[wt/math,text=wtblack,anchor=north] at (6.700,-4.570) {$\displaystyle a_{k-1}$};
\node[wt/math,text=wtblack,anchor=north] at (8.400,-4.570) {$\displaystyle a_k$};
\node[wt/math,text=wtblack,anchor=north] at (14.500,-4.570) {$\displaystyle a_m$};
\node[wt/math,text=wtblack,anchor=east] at (0.330,-3.480) {$\displaystyle b$};
\node[wt/math,text=wtblack,anchor=east] at (0.330,-5.120) {$\displaystyle c$};
\node[wt/math,text=wtblue,anchor=center] at (8.000,-2.770) {$\displaystyle F=\{b,c\},\quad s=2$};
\node[wt/math,text=wtblue,anchor=south] at (4.600,-3.830) {$\displaystyle S$};
\node[wt/math,text=wtblue,anchor=north] at (4.450,-5.510) {$\displaystyle k$};
\node[wt/math,text=wtblue,anchor=north] at (11.450,-5.510) {$\displaystyle m+1-k$};
\node[wt/math,text=wtblue,anchor=center] at (8.000,-6.220) {$\displaystyle S=\{a_0,\ldots,a_{k-1}\}\cup F,\quad |S|=k+s$};
\end{tikzpicture}
}
\caption{Connected averaging sets containing one or two leaves.}
\label{fig:selected-leaves}
\end{figure}
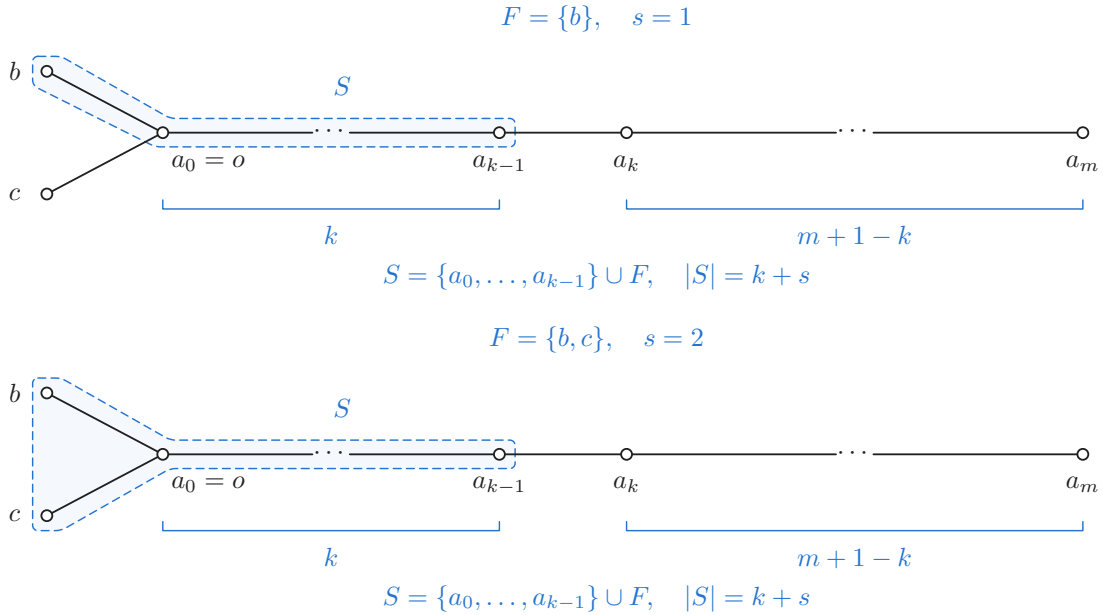

Since
\[
T_{m+1-k}
=\{a_k,a_{k+1},\ldots,a_m\}
\subseteq V(S(m,1,1))\setminus S,
\]
and all weights are nonnegative,
\[
\begin{aligned}
\sum_{v\in S}\omega(v)
&=1-\sum_{v\in V(S(m,1,1))\setminus S}\omega(v)\\
&\leq 1-\sum_{v\in T_{m+1-k}}\omega(v)\\
&=1-Y_{m+1-k}.
\end{aligned}
\]
Let $\mu$ be the common weight on $S$ after the operation. As $|S|=k+s$,
\[
\mu
=\frac{1}{k+s}\sum_{v\in S}\omega(v)
\leq \frac{1-Y_{m+1-k}}{k+s}.
\]

If $0\le t\le m+1-k$, then
\[
T_t\subseteq T_{m+1-k}
=\{a_k,a_{k+1},\ldots,a_m\},
\]
so $T_t\cap S=\varnothing$. No weight on $T_t$ changes, and hence
\[
Y_t'=Y_t\leq \frac{t}{m+1}.
\]

If $m+1-k<t\le m+1$, define
\[
h=t-(m+1-k).
\]
Thus $h$ counts the vertices of $T_t\cap S$. Since
\[
m+1-k<t\leq m+1
\]
we have
\[
1\leq h\leq k.
\]
Moreover,
\[
T_t\setminus S
=\{a_k,a_{k+1},\ldots,a_m\}
=T_{m+1-k},
\]
and
\[
T_t\cap S
=\{a_{k-h},a_{k-h+1},\ldots,a_{k-1}\}.
\]
See Figure~\ref{fig:intersection-h} for the intersection and terminal
segments involved.

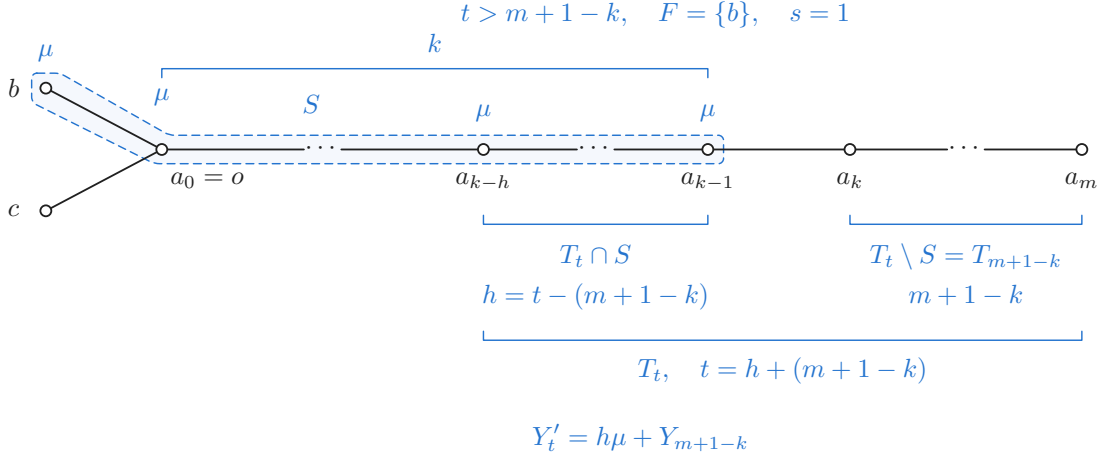
\begin{figure}[!htbp]
\centering
\resizebox{.96\linewidth}{!}{
\begin{tikzpicture}[wt/base]
\path (0,-4.5) rectangle (15.4,2.0);
\draw[wt/selection] (0.460,1.020) -- (0.840,1.020) -- (2.300,0.190) -- (9.710,0.190) -- (9.710,-0.190) -- (2.100,-0.190) -- (0.460,0.630) -- cycle;
\draw[wt/edge] (0.650,0.820) -- (2.200,0.000);
\draw[wt/edge] (0.650,-0.820) -- (2.200,0.000);
\draw[wt/edge] (2.200,0.000) -- (4.100,0.000);
\draw[wt/edge] (4.600,0.000) -- (6.500,0.000);
\draw[wt/edge] (6.500,0.000) -- (7.750,0.000);
\draw[wt/edge] (8.250,0.000) -- (9.500,0.000);
\draw[wt/edge] (9.500,0.000) -- (11.400,0.000);
\draw[wt/edge] (11.400,0.000) -- (12.700,0.000);
\draw[wt/edge] (13.200,0.000) -- (14.500,0.000);
\draw[wt/bracket,draw=wtblue] (2.200,0.960) -- (2.200,1.080) -- (9.500,1.080) -- (9.500,0.960);
\draw[wt/bracket,draw=wtblue] (6.500,-0.880) -- (6.500,-1.000) -- (9.500,-1.000) -- (9.500,-0.880);
\draw[wt/bracket,draw=wtblue] (11.400,-0.880) -- (11.400,-1.000) -- (14.500,-1.000) -- (14.500,-0.880);
\draw[wt/bracket,draw=wtblue] (6.500,-2.430) -- (6.500,-2.550) -- (14.500,-2.550) -- (14.500,-2.430);
\draw[draw=wtblack,fill=white,line width=.70pt] (2.200,0.000) circle[radius=.073cm];
\draw[draw=wtblack,fill=white,line width=.70pt] (6.500,0.000) circle[radius=.073cm];
\draw[draw=wtblack,fill=white,line width=.70pt] (9.500,0.000) circle[radius=.073cm];
\draw[draw=wtblack,fill=white,line width=.70pt] (11.400,0.000) circle[radius=.073cm];
\draw[draw=wtblack,fill=white,line width=.70pt] (14.500,0.000) circle[radius=.073cm];
\draw[wt/edge,fill=white] (0.650,0.820) circle[radius=.073cm];
\draw[wt/edge,fill=white] (0.650,-0.820) circle[radius=.073cm];
\node[wt/math,text=wtblack,anchor=center] at (4.350,0.018) {$\displaystyle \cdots$};
\node[wt/math,text=wtblack,anchor=center] at (8.000,0.018) {$\displaystyle \cdots$};
\node[wt/math,text=wtblack,anchor=center] at (12.950,0.018) {$\displaystyle \cdots$};
\node[wt/math,text=wtblue,anchor=south] at (2.200,0.560) {$\displaystyle \mu$};
\node[wt/math,text=wtblack,anchor=north west] at (2.280,-0.270) {$\displaystyle a_0=o$};
\node[wt/math,text=wtblack,anchor=north] at (6.500,-0.270) {$\displaystyle a_{k-h}$};
\node[wt/math,text=wtblue,anchor=south] at (6.500,0.340) {$\displaystyle \mu$};
\node[wt/math,text=wtblack,anchor=north] at (9.500,-0.270) {$\displaystyle a_{k-1}$};
\node[wt/math,text=wtblue,anchor=south] at (9.500,0.340) {$\displaystyle \mu$};
\node[wt/math,text=wtblack,anchor=north] at (11.400,-0.270) {$\displaystyle a_k$};
\node[wt/math,text=wtblack,anchor=north] at (14.500,-0.270) {$\displaystyle a_m$};
\node[wt/math,text=wtblack,anchor=east] at (0.330,0.820) {$\displaystyle b$};
\node[wt/math,text=wtblue,anchor=south] at (0.650,1.140) {$\displaystyle \mu$};
\node[wt/math,text=wtblack,anchor=east] at (0.330,-0.820) {$\displaystyle c$};
\node[wt/math,text=wtblue,anchor=center] at (8.800,1.800) {$\displaystyle t>m+1-k,\quad F=\{b\},\quad s=1$};
\node[wt/math,text=wtblue,anchor=south] at (5.850,1.270) {$\displaystyle k$};
\node[wt/math,text=wtblue,anchor=south] at (4.200,0.450) {$\displaystyle S$};
\node[wt/math,text=wtblue,anchor=north] at (8.000,-1.190) {$\displaystyle \begin{gathered}T_t\cap S\\ h=t-(m+1-k)\end{gathered}$};
\node[wt/math,text=wtblue,anchor=north] at (12.950,-1.190) {$\displaystyle \begin{gathered}T_t\setminus S=T_{m+1-k}\\ m+1-k\end{gathered}$};
\node[wt/math,text=wtblue,anchor=north] at (10.500,-2.740) {$\displaystyle T_t,\quad t=h+(m+1-k)$};
\node[wt/math,text=wtblue,anchor=center] at (8.600,-3.850) {$\displaystyle Y_t'=h\mu+Y_{m+1-k}$};
\end{tikzpicture}
}
\caption{The intersection count $h$ and the decomposition of $T_t$.}
\label{fig:intersection-h}
\end{figure}

The operation leaves the weights on $T_{m+1-k}$ unchanged and assigns
weight $\mu$ to each of the $h$ vertices of $T_t\cap S$. Thus
\begin{equation}\label{eq:root-leaf-strict}
\begin{aligned}
Y_t'
&=Y_{m+1-k}+h\mu\\
&\leq
Y_{m+1-k}
+\frac{h}{k+s}\bigl(1-Y_{m+1-k}\bigr)\\
&=
\left(1-\frac{h}{k+s}\right)Y_{m+1-k}
+\frac{h}{k+s}\\
&\leq
\left(1-\frac{h}{k+s}\right)\frac{m+1-k}{m+1}
+\frac{h}{k+s}\\
&=\frac{t}{m+1}
-\frac{hs}{(k+s)(m+1)}\\
&<\frac{t}{m+1}.
\end{aligned}
\end{equation}
The fourth line uses \eqref{eq:root-tail-bound}, and
\[
1-\frac{h}{k+s}>0
\]
because $h\le k<k+s$. This proves \eqref{eq:root-tail-bound} for $\omega'$.

To prove \eqref{eq:root-tail-rigidity}, suppose that $1\le t\le m$ and
\[
Y_t'=\frac{t}{m+1}.
\]
By \eqref{eq:root-leaf-strict}, we cannot have $t>m+1-k$. Hence
\[
t\leq m+1-k.
\]
In this case $T_t\cap S=\varnothing$, so no weight on $T_t$ changes.
By \eqref{eq:root-tail-rigidity}, before the operation,
\[
\omega(a_{m-t+1})=\cdots=\omega(a_m)=\frac{1}{m+1},
\]
and therefore, after the operation,
\[
\omega'(a_{m-t+1})=\cdots=\omega'(a_m)=\frac{1}{m+1}.
\]
Thus $\omega'$ also satisfies both \eqref{eq:root-tail-bound} and \eqref{eq:root-tail-rigidity} when
$S\cap\{b,c\}\ne\varnothing$.

The induction is complete. Thus \eqref{eq:root-tail-bound} and \eqref{eq:root-tail-rigidity} hold after any finite sequence
of averaging operations. Taking $t=1$ in \eqref{eq:root-tail-bound} gives
\[
\omega(a_m)=Y_1\leq \frac{1}{m+1}.
\]

It remains to determine when equality holds. Suppose that, after finitely
many averaging operations,
\[
\omega(a_m)=\frac{1}{m+1}.
\]
Since the initial weight at $a_m$ is zero, some effective operation must
contain $a_m$. Choose the last such operation, denote its selected set by
$S$, and write $\omega$ and $\omega'$ for the weights immediately before and after it,
respectively.

No later effective operation contains $a_m$. Operations avoiding $a_m$
and ineffective operations containing it leave its weight unchanged.
Consequently, immediately after the chosen operation,
\[
\omega'(a_m)=\frac{1}{m+1}.
\]

If
\[
S\cap\{b,c\}\neq\varnothing,
\]
then connectivity and $a_m\in S$ imply
\[
\{a_0,a_1,\ldots,a_m\}\subseteq S.
\]
Thus
\[
|S|\geq m+2.
\]
The weights before the operation are nonnegative and have total sum one,
so the common weight on $S$ after the operation satisfies
\[
\omega'(a_m)
=\frac{1}{|S|}\sum_{v\in S}\omega(v)
\leq \frac{1}{m+2}
<\frac{1}{m+1}.
\]
This contradicts $\omega'(a_m)=1/(m+1)$. Therefore,
\[
S\cap\{b,c\}=\varnothing.
\]

Since $S$ is connected and contains $a_m$, there is an integer
$0\le r\le m$ such that
\[
S=\{a_r,a_{r+1},\ldots,a_m\}.
\]
If $r\ge1$, put
\[
t=m-r+1.
\]
Thus $t=|S|$ and $1\le t\le m$. After the operation, every vertex of $S$
has weight $\omega'(a_m)=1/(m+1)$. Since the operation preserves the total
weight on $S$, before the operation we have
\[
Y_t=\sum_{i=r}^{m}\omega(a_i)=\frac{t}{m+1}.
\]
By \eqref{eq:root-tail-rigidity},
\[
\omega(a_r)=\omega(a_{r+1})=\cdots=\omega(a_m)=\frac{1}{m+1}.
\]
The weights on $S$ are therefore already equal before the operation,
contrary to its effectiveness. Hence
\[
r=0,
\qquad
S=\{a_0,a_1,\ldots,a_m\}.
\]

\begin{figure}[!htbp]
\centering
\resizebox{.96\linewidth}{!}{
\begin{tikzpicture}[wt/base]
\path (0,-2.1) rectangle (15.4,1.9);
\draw[wt/selection] (1.990,-0.190) rectangle (14.710,0.190);
\draw[wt/edge] (0.650,0.820) -- (2.200,0.000);
\draw[wt/edge] (0.650,-0.820) -- (2.200,0.000);
\draw[wt/edge] (2.200,0.000) -- (6.900,0.000);
\draw[wt/edge] (6.900,0.000) -- (10.450,0.000);
\draw[wt/edge] (10.950,0.000) -- (14.500,0.000);
\draw[draw=wtblack,fill=white,line width=.70pt] (2.200,0.000) circle[radius=.073cm];
\draw[draw=wtblack,fill=white,line width=.70pt] (6.900,0.000) circle[radius=.073cm];
\draw[draw=wtblack,fill=white,line width=.70pt] (14.500,0.000) circle[radius=.073cm];
\draw[wt/edge,fill=white] (0.650,0.820) circle[radius=.073cm];
\draw[wt/edge,fill=white] (0.650,-0.820) circle[radius=.073cm];
\node[wt/math,text=wtblack,anchor=center] at (10.700,0.018) {$\displaystyle \cdots$};
\node[wt/math,text=wtblue,anchor=south] at (2.200,0.340) {$\displaystyle \frac1{m+1}$};
\node[wt/math,text=wtblack,anchor=north west] at (2.280,-0.270) {$\displaystyle a_0=o$};
\node[wt/math,text=wtblack,anchor=north] at (6.900,-0.270) {$\displaystyle a_1$};
\node[wt/math,text=wtblue,anchor=south] at (6.900,0.340) {$\displaystyle \frac1{m+1}$};
\node[wt/math,text=wtblack,anchor=north] at (14.500,-0.270) {$\displaystyle a_m$};
\node[wt/math,text=wtblue,anchor=south] at (14.500,0.340) {$\displaystyle \frac1{m+1}$};
\node[wt/math,text=wtblack,anchor=east] at (0.330,0.820) {$\displaystyle b$};
\node[wt/math,text=wtblue,anchor=south] at (0.650,1.140) {$\displaystyle 0$};
\node[wt/math,text=wtblack,anchor=east] at (0.330,-0.820) {$\displaystyle c$};
\node[wt/math,text=wtblue,anchor=north] at (0.650,-1.140) {$\displaystyle 0$};
\node[wt/math,text=wtblue,anchor=center] at (8.600,1.550) {$\displaystyle \omega'(a_m)=\frac1{m+1}$};
\node[wt/math,text=wtblue,anchor=center] at (8.600,-1.430) {$\displaystyle S=\{a_0,a_1,\ldots,a_m\},\quad |S|=m+1$};
\end{tikzpicture}
}
\caption{The equality configuration for the root-source barrier.}
\label{fig:root-equality}
\end{figure}
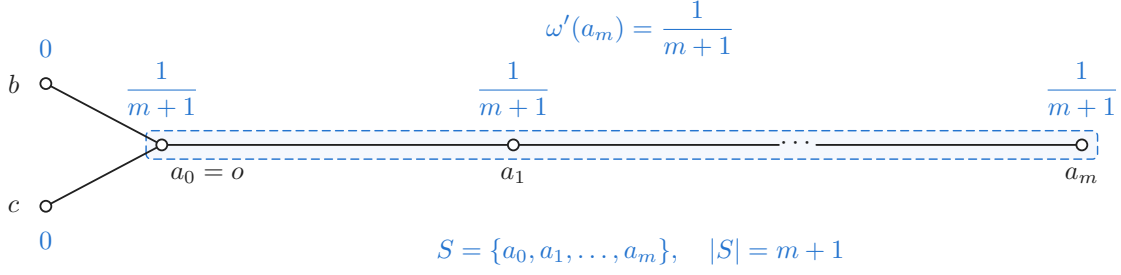

Immediately after the operation, all $m+1$ vertices of $S$ have weight
$1/(m+1)$. Preservation of the total weight on $S$ gives, before the operation,
\[
\sum_{i=0}^{m}\omega(a_i)
=(m+1)\frac{1}{m+1}
=1.
\]
Since the total weight is one and $\omega(b),\omega(c)\ge0$,
\[
\omega(b)=\omega(c)=0.
\]
Neither $b$ nor $c$ is selected, so
\[
\omega'(b)=\omega'(c)=0.
\]
Thus, immediately after the operation,
\[
\omega'(a_0)=\omega'(a_1)=\cdots=\omega'(a_m)=\frac{1}{m+1},
\qquad
\omega'(b)=\omega'(c)=0.
\]
The equality configuration is shown in Figure~\ref{fig:root-equality}.

\end{proof}

\FloatBarrier
\begin{lemma}[Barrier for two leaf sources]
Suppose that the initial weight function $\omega$ on $S(m,1,1)$ satisfies
\[
\omega(a_i)=0 \quad (0\le i\le m),
\qquad
\omega(b)=\omega(c)=1.
\]
Then, after any finite sequence of averaging operations,
\[
\omega(a_m)\le \frac{2m+3}{(m+2)^2}.
\]
If equality holds, then the last effective operation containing $a_m$ must
average one of the sets
\[
\{a_0,a_1,\ldots,a_m,b\}
\quad\text{ or }\quad
\{a_0,a_1,\ldots,a_m,c\}.
\]
In the first case, immediately after that operation,
\[
\omega(a_0)=\omega(a_1)=\cdots=\omega(a_m)=\omega(b)
=\frac{2m+3}{(m+2)^2},
\qquad
\omega(c)=\frac{1}{m+2}.
\]
In the second case, immediately after that operation,
\[
\omega(a_0)=\omega(a_1)=\cdots=\omega(a_m)=\omega(c)
=\frac{2m+3}{(m+2)^2},
\qquad
\omega(b)=\frac{1}{m+2}.
\]
Moreover, every sequence attaining equality contains at least two
effective operations.
\end{lemma}

\begin{proof}
As in the proof of Lemma~1, put
\[
T_t=\{a_{m-t+1},a_{m-t+2},\ldots,a_m\},
\qquad 1\le t\le m+1,
\]
and $T_0=\varnothing$. Thus $T_t$ consists of the last $t$ vertices of
$\langle a_0,a_1,\ldots,a_m\rangle$, counting from $a_m$ towards $a_0$.
In particular,
\[
Y_t=\sum_{v\in T_t}\omega(v).
\]

Averaging preserves nonnegativity and the sum of all vertex weights.
Hence, throughout the sequence,
\[
\omega(v)\ge 0
\quad\bigl(v\in V(S(m,1,1))\bigr),
\qquad
\sum_{v\in V(S(m,1,1))}\omega(v)=2.
\]

We prove by induction on the number of operations that
\begin{equation}\label{eq:two-leaf-tail-bound}
Y_t\le \frac{t(2m+3)}{(m+2)^2},
\qquad 0\le t\le m+1,
\end{equation}
and
\begin{equation}\label{eq:two-leaf-prefix-bound}
\omega(b)+X_k\ge \frac{k+1}{m+2},
\qquad
\omega(c)+X_k\ge \frac{k+1}{m+2},
\qquad 0\le k\le m+1,
\end{equation}
and, simultaneously, that for $1\le t\le m+1$,
\begin{equation}\label{eq:two-leaf-tail-rigidity}
Y_t=\frac{t(2m+3)}{(m+2)^2}
\quad\Longrightarrow\quad
\omega(a_{m-t+1})=\cdots=\omega(a_m)
=\frac{2m+3}{(m+2)^2}.
\end{equation}
Taking $t=1$ in \eqref{eq:two-leaf-tail-bound} gives
\[
\omega(a_m)=Y_1\le \frac{2m+3}{(m+2)^2},
\]
which is the required bound. We use \eqref{eq:two-leaf-prefix-bound} to show that averaging preserves
\eqref{eq:two-leaf-tail-bound}, and \eqref{eq:two-leaf-tail-rigidity} to identify the last effective operation containing $a_m$
when equality holds.

Initially, every path vertex has weight zero, so
\[
X_k=0 \quad (0\le k\le m+1),
\qquad
Y_t=0 \quad (0\le t\le m+1),
\]
and \eqref{eq:two-leaf-tail-bound} holds. Also,
\[
\omega(b)=\omega(c)=1\ge \frac{k+1}{m+2},
\qquad 0\le k\le m+1,
\]
so \eqref{eq:two-leaf-prefix-bound} holds. For $1\le t\le m+1$, the equality in the antecedent of \eqref{eq:two-leaf-tail-rigidity}
does not hold initially. Thus \eqref{eq:two-leaf-tail-rigidity} holds as well.

Suppose that the current weight function $\omega$ satisfies \eqref{eq:two-leaf-tail-bound}--\eqref{eq:two-leaf-tail-rigidity}.
Let $S$ be the connected set selected by the next operation, let $\omega'$ be
the resulting weight function, and write
\[
X_k'=\sum_{i=0}^{k-1}\omega'(a_i),
\qquad
Y_t'=\sum_{v\in T_t}\omega'(v),
\]
with $X_0'=Y_0'=0$. Thus $X_k'$ is the sum of the new weights on the
first $k$ path vertices starting at $a_0$, and $Y_t'$ is the sum of the
new weights on $T_t$.

If the operation is ineffective, then $\omega'=\omega$ and the conclusion follows.
We may therefore assume that it is effective.

First suppose that
\[
S\cap\{b,c\}=\varnothing.
\]
Since $S$ is connected, there are integers $0\le r\le j\le m$ such that
\[
S=\{a_r,a_{r+1},\ldots,a_j\}.
\]
Apply the interpolation identity \eqref{eq:path-interpolation} and its equality analysis from the
proof of Lemma~1, using the induction hypotheses \eqref{eq:two-leaf-tail-bound} and \eqref{eq:two-leaf-tail-rigidity}.
It follows that $\omega'$ still satisfies \eqref{eq:two-leaf-tail-bound} and \eqref{eq:two-leaf-tail-rigidity}.

We verify \eqref{eq:two-leaf-prefix-bound}. Fix $u\in\{b,c\}$ and an integer $0\le p\le m+1$.
Since neither leaf is selected, $\omega'(u)=\omega(u)$.

If $p\le r$, then $\{a_0,\ldots,a_{p-1}\}$ is disjoint from $S$.
If $p\ge j+1$, then this set contains $S$. In either case,
$X_p'=X_p$, so the induction hypothesis \eqref{eq:two-leaf-prefix-bound} gives
\[
\omega'(u)+X_p'=\omega(u)+X_p\ge\frac{p+1}{m+2}.
\]

If $r<p<j+1$, the operation preserves the total weight on the path,
so $X_{m+1}'=X_{m+1}$. Using
\[
X_p'=X_{m+1}-Y_{m+1-p}',
\]
apply the interpolation identity \eqref{eq:path-interpolation} to $Y_{m+1-p}'$ and then use the
induction hypothesis \eqref{eq:two-leaf-prefix-bound} to obtain
\[
\begin{aligned}
\omega'(u)+X_p'
&=\frac{j-p+1}{j-r+1}\bigl(\omega(u)+X_r\bigr)
 +\frac{p-r}{j-r+1}\bigl(\omega(u)+X_{j+1}\bigr)\\
&\ge
\frac{j-p+1}{j-r+1}\frac{r+1}{m+2}
+
\frac{p-r}{j-r+1}\frac{j+2}{m+2}\\
&=\frac{p+1}{m+2}.
\end{aligned}
\]
Thus $\omega'$ satisfies \eqref{eq:two-leaf-tail-bound}--\eqref{eq:two-leaf-tail-rigidity} when $S\cap\{b,c\}=\varnothing$.

Next suppose that
\[
S\cap\{b,c\}\ne\varnothing.
\]
Define
\[
F=S\cap\{b,c\},
\qquad
s=|F|,
\qquad
k=\bigl|S\cap\{a_0,a_1,\ldots,a_m\}\bigr|.
\]
Thus $F$ is the set of selected vertices among $b,c$, $s$ is their number,
and $k$ is the number of selected vertices on the path
$\langle a_0,a_1,\ldots,a_m\rangle$.

Since the operation is effective and $S$ is connected,
\[
s\in\{1,2\},
\qquad
1\le k\le m+1,
\qquad
S=\{a_0,a_1,\ldots,a_{k-1}\}\cup F.
\]
Let $\mu$ denote the common weight on $S$ after the operation.

We first use \eqref{eq:two-leaf-prefix-bound} to prove that
\[
\mu\ge \frac{1}{m+2}.
\]
If $s=1$, write $F=\{u\}$. Then
\[
\mu=\frac{\omega(u)+X_k}{k+1}
\ge \frac{1}{m+2}.
\]
If $s=2$, apply \eqref{eq:two-leaf-prefix-bound} with indices $k$ and $0$ to obtain
\[
\mu
=\frac{\omega(b)+X_k+\omega(c)}{k+2}
\ge
\frac{\dfrac{k+1}{m+2}+\dfrac{1}{m+2}}{k+2}
=\frac{1}{m+2}.
\]

We now verify \eqref{eq:two-leaf-prefix-bound}. Fix $u\in\{b,c\}$ and an integer $0\le p\le m+1$.
If $p\le k$, then
\[
X_p'=p\mu.
\]
When $u\in F$,
\[
\omega'(u)+X_p'=(p+1)\mu
\ge \frac{p+1}{m+2}.
\]
When $u\notin F$, apply \eqref{eq:two-leaf-prefix-bound} with index $0$ to obtain
\[
\omega'(u)+X_p'
=\omega(u)+p\mu
\ge \frac{1}{m+2}+\frac{p}{m+2}
=\frac{p+1}{m+2}.
\]

Now suppose that $p>k$, and put
\[
q=p-k,
\qquad
t=m+1-p.
\]
Here $q$ counts the vertices of $\{a_k,a_{k+1},\ldots,a_{p-1}\}$,
and $t$ counts the path vertices outside $\{a_0,a_1,\ldots,a_{p-1}\}$.
Thus
\[
q\ge1,\qquad t\ge0,\qquad m+1=k+q+t.
\]

If $s=1$ and $u\in F$, then
\[
S\subseteq\{u,a_0,a_1,\ldots,a_{p-1}\}.
\]
The operation preserves the total weight on the set on the right, so
\[
\omega'(u)+X_p'=\omega(u)+X_p
\ge \frac{p+1}{m+2}.
\]

It remains to consider $s=2$, or $s=1$ and $u\notin F$.
Let $\bar u$ be the unique vertex of $\{b,c\}\setminus\{u\}$.
In either case, $\bar u\in F$, and hence
\[
\omega'(\bar u)=\mu.
\]
The complement of $\{u,a_0,a_1,\ldots,a_{p-1}\}$ in the full vertex set is
\[
\{\bar u\}\cup T_t.
\]
Since $p>k$, we have $T_t\cap S=\varnothing$ and thus $Y_t'=Y_t$.
The total weight is two, so
\[
2-\bigl(\omega'(u)+X_p'\bigr)=\mu+Y_t.
\]
The counts and complementary sets are illustrated in
Figure~\ref{fig:pqt-complement}.

\begin{figure}[!htbp]
\centering
\resizebox{.96\linewidth}{!}{
\begin{tikzpicture}[wt/base]
\path (0,-3.4) rectangle (15.4,2.15);
\draw[wt/selection] (0.460,1.020) -- (0.840,1.020) -- (2.300,0.190) -- (5.310,0.190) -- (5.310,-0.190) -- (2.100,-0.190) -- (0.460,0.630) -- cycle;
\draw[wt/edge] (0.650,0.820) -- (2.200,0.000);
\draw[wt/edge] (0.650,-0.820) -- (2.200,0.000);
\draw[wt/edge] (2.200,0.000) -- (3.400,0.000);
\draw[wt/edge] (3.900,0.000) -- (5.100,0.000);
\draw[wt/edge] (5.100,0.000) -- (6.700,0.000);
\draw[wt/edge] (6.700,0.000) -- (7.850,0.000);
\draw[wt/edge] (8.350,0.000) -- (9.500,0.000);
\draw[wt/edge] (9.500,0.000) -- (11.100,0.000);
\draw[wt/edge] (11.100,0.000) -- (12.550,0.000);
\draw[wt/edge] (13.050,0.000) -- (14.500,0.000);
\draw[wt/bracket,draw=wtblue] (2.200,0.950) -- (2.200,1.070) -- (5.100,1.070) -- (5.100,0.950);
\draw[wt/bracket,draw=wtblue] (6.700,0.950) -- (6.700,1.070) -- (9.500,1.070) -- (9.500,0.950);
\draw[wt/bracket,draw=wtblue] (11.100,0.950) -- (11.100,1.070) -- (14.500,1.070) -- (14.500,0.950);
\draw[wt/bracket,draw=wtblue] (2.200,-0.880) -- (2.200,-1.000) -- (9.500,-1.000) -- (9.500,-0.880);
\draw[wt/bracket,draw=wtblue] (11.100,-0.880) -- (11.100,-1.000) -- (14.500,-1.000) -- (14.500,-0.880);
\draw[draw=wtblack,fill=white,line width=.70pt] (2.200,0.000) circle[radius=.073cm];
\draw[draw=wtblack,fill=white,line width=.70pt] (5.100,0.000) circle[radius=.073cm];
\draw[draw=wtblack,fill=white,line width=.70pt] (6.700,0.000) circle[radius=.073cm];
\draw[draw=wtblack,fill=white,line width=.70pt] (9.500,0.000) circle[radius=.073cm];
\draw[draw=wtblack,fill=white,line width=.70pt] (11.100,0.000) circle[radius=.073cm];
\draw[draw=wtblack,fill=white,line width=.70pt] (14.500,0.000) circle[radius=.073cm];
\draw[wt/edge,fill=white] (0.650,0.820) circle[radius=.073cm];
\draw[wt/edge,fill=white] (0.650,-0.820) circle[radius=.073cm];
\node[wt/math,text=wtblack,anchor=center] at (3.650,0.018) {$\displaystyle \cdots$};
\node[wt/math,text=wtblack,anchor=center] at (8.100,0.018) {$\displaystyle \cdots$};
\node[wt/math,text=wtblack,anchor=center] at (12.800,0.018) {$\displaystyle \cdots$};
\node[wt/math,text=wtblue,anchor=south] at (2.200,0.560) {$\displaystyle \mu$};
\node[wt/math,text=wtblack,anchor=north west] at (2.280,-0.270) {$\displaystyle a_0=o$};
\node[wt/math,text=wtblack,anchor=north] at (5.100,-0.270) {$\displaystyle a_{k-1}$};
\node[wt/math,text=wtblue,anchor=south] at (5.100,0.340) {$\displaystyle \mu$};
\node[wt/math,text=wtblack,anchor=north] at (6.700,-0.270) {$\displaystyle a_k$};
\node[wt/math,text=wtblack,anchor=north] at (9.500,-0.270) {$\displaystyle a_{p-1}$};
\node[wt/math,text=wtblack,anchor=north] at (11.100,-0.270) {$\displaystyle a_p$};
\node[wt/math,text=wtblack,anchor=north] at (14.500,-0.270) {$\displaystyle a_m$};
\node[wt/math,text=wtblack,anchor=east] at (0.330,0.820) {$\displaystyle \bar u$};
\node[wt/math,text=wtblue,anchor=south] at (0.650,1.140) {$\displaystyle \mu$};
\node[wt/math,text=wtblack,anchor=east] at (0.330,-0.820) {$\displaystyle u$};
\node[wt/math,text=wtblue,anchor=north] at (0.650,-1.140) {$\displaystyle \omega'(u)$};
\node[wt/math,text=wtblue,anchor=center] at (8.800,1.910) {$\displaystyle p>k,\quad F=\{\bar u\},\quad s=1$};
\node[wt/math,text=wtblue,anchor=south] at (3.650,1.260) {$\displaystyle k$};
\node[wt/math,text=wtblue,anchor=south] at (8.100,1.260) {$\displaystyle q=p-k$};
\node[wt/math,text=wtblue,anchor=south] at (12.800,1.260) {$\displaystyle t=m+1-p$};
\node[wt/math,text=wtblue,anchor=south] at (3.900,0.450) {$\displaystyle S$};
\node[wt/math,text=wtblue,anchor=north] at (5.850,-1.190) {$\displaystyle X_p',\quad p=k+q$};
\node[wt/math,text=wtblue,anchor=north] at (12.800,-1.190) {$\displaystyle \begin{gathered}T_t\\Y_t'=Y_t\end{gathered}$};
\node[wt/math,text=wtblue,anchor=center] at (8.100,-2.120) {$\displaystyle V(S(m,1,1))\setminus\{u,a_0,\ldots,a_{p-1}\}=\{\bar u\}\cup T_t$};
\node[wt/math,text=wtblue,anchor=center] at (8.100,-2.970) {$\displaystyle 2-\bigl(\omega'(u)+X_p'\bigr)=\mu+Y_t$};
\end{tikzpicture}
}
\caption{The counts $p,k,q,t$ and the corresponding complementary sets.}
\label{fig:pqt-complement}
\end{figure}
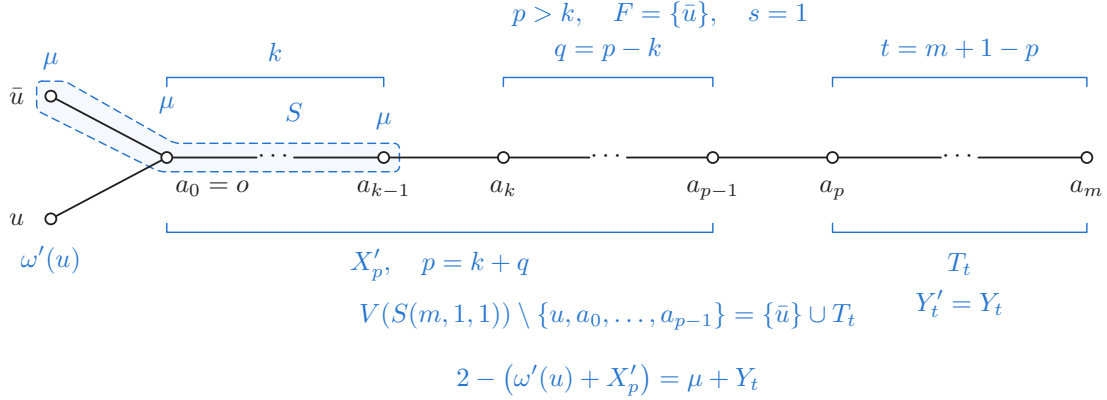

If $s=1$ and $u\notin F$, then
\[
\mu=\frac{2-\omega(u)-Y_{m+1-k}}{k+1}.
\]
By \eqref{eq:two-leaf-prefix-bound}, $\omega(u)\ge1/(m+2)$. Also, $T_t\subseteq T_{m+1-k}$ and all
weights are nonnegative, so
\[
Y_{m+1-k}\ge Y_t.
\]
It follows from \eqref{eq:two-leaf-tail-bound} that
\[
\begin{aligned}
\mu+Y_t
&\le \frac{2-\dfrac{1}{m+2}+kY_t}{k+1}\\
&\le
\frac{2-\dfrac{1}{m+2}
      +\dfrac{kt(2m+3)}{(m+2)^2}}{k+1}\\
&\le
\frac{2-\dfrac{1}{m+2}+\dfrac{2kt}{m+2}}{k+1}\\
&=
2-\frac{p+1}{m+2}
-\frac{k^2+(k-1)q}{(k+1)(m+2)}\\
&<
2-\frac{p+1}{m+2}.
\end{aligned}
\]
The third line uses
\[
\frac{2m+3}{(m+2)^2}
=\frac{2}{m+2}-\frac{1}{(m+2)^2}
<\frac{2}{m+2},
\]
the fourth uses $p=k+q$ and $m+1=k+q+t$, and the last inequality follows
from $k\ge1$ and $q\ge1$. Hence
\[
\omega'(u)+X_p'>\frac{p+1}{m+2}.
\]

If $s=2$, then
\[
\mu=\frac{2-Y_{m+1-k}}{k+2}.
\]
Again using $Y_{m+1-k}\ge Y_t$ and \eqref{eq:two-leaf-tail-bound}, we obtain
\[
\begin{aligned}
\mu+Y_t
&\le \frac{2+(k+1)Y_t}{k+2}\\
&\le
\frac{2+\dfrac{(k+1)t(2m+3)}{(m+2)^2}}{k+2}\\
&\le
\frac{2+\dfrac{2(k+1)t}{m+2}}{k+2}\\
&=
2-\frac{p+1}{m+2}
-\frac{k(k+q+1)}{(k+2)(m+2)}\\
&<
2-\frac{p+1}{m+2}.
\end{aligned}
\]
Therefore,
\[
\omega'(u)+X_p'>\frac{p+1}{m+2}.
\]
This proves \eqref{eq:two-leaf-prefix-bound} for $\omega'$.

We next verify \eqref{eq:two-leaf-tail-bound} and \eqref{eq:two-leaf-tail-rigidity}. Define
\[
Z=\sum_{u\in\{b,c\}\setminus F}\omega(u).
\]
Thus $Z$ is the total weight on the unselected vertices among $b,c$.
If $s=1$, apply \eqref{eq:two-leaf-prefix-bound} with index $0$ to obtain
\[
Z\ge\frac{1}{m+2}.
\]
If $s=2$, then $Z=0$.

The vertices outside $S$ consist precisely of $T_{m+1-k}$ and
$\{b,c\}\setminus F$. Hence
\[
\mu=\frac{2-Y_{m+1-k}-Z}{k+s}.
\]
These two contributions are shown in Figure~\ref{fig:outside-mass}.

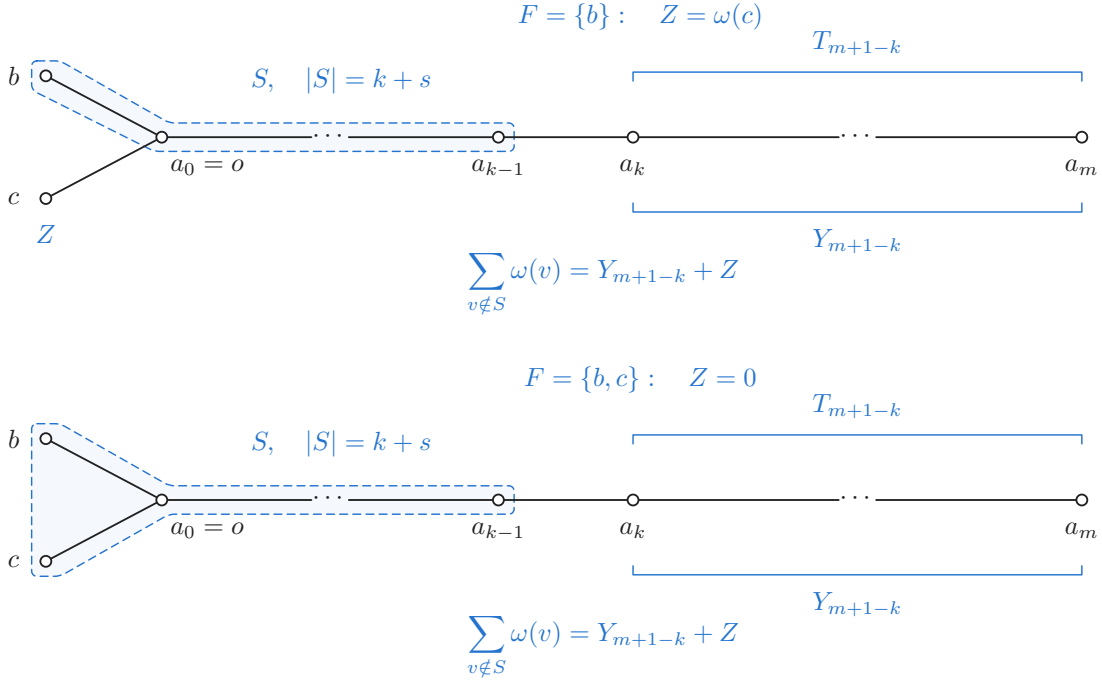
\begin{figure}[!htbp]
\centering
\resizebox{.96\linewidth}{!}{
\begin{tikzpicture}[wt/base]
\path (0,-7.35) rectangle (15.4,1.9);
\draw[wt/selection] (0.460,1.020) -- (0.840,1.020) -- (2.300,0.190) -- (6.910,0.190) -- (6.910,-0.190) -- (2.100,-0.190) -- (0.460,0.630) -- cycle;
\draw[wt/selection] (0.460,-3.830) -- (0.840,-3.830) -- (2.300,-4.660) -- (6.910,-4.660) -- (6.910,-5.040) -- (2.300,-5.040) -- (0.840,-5.870) -- (0.460,-5.870) -- cycle;
\draw[wt/edge] (0.650,0.820) -- (2.200,0.000);
\draw[wt/edge] (0.650,-0.820) -- (2.200,0.000);
\draw[wt/edge] (2.200,0.000) -- (4.200,0.000);
\draw[wt/edge] (4.700,0.000) -- (6.700,0.000);
\draw[wt/edge] (6.700,0.000) -- (8.500,0.000);
\draw[wt/edge] (8.500,0.000) -- (11.250,0.000);
\draw[wt/edge] (11.750,0.000) -- (14.500,0.000);
\draw[wt/bracket,draw=wtblue] (8.500,0.750) -- (8.500,0.870) -- (14.500,0.870) -- (14.500,0.750);
\draw[wt/bracket,draw=wtblue] (8.500,-0.880) -- (8.500,-1.000) -- (14.500,-1.000) -- (14.500,-0.880);
\draw[wt/edge] (0.650,-4.030) -- (2.200,-4.850);
\draw[wt/edge] (0.650,-5.670) -- (2.200,-4.850);
\draw[wt/edge] (2.200,-4.850) -- (4.200,-4.850);
\draw[wt/edge] (4.700,-4.850) -- (6.700,-4.850);
\draw[wt/edge] (6.700,-4.850) -- (8.500,-4.850);
\draw[wt/edge] (8.500,-4.850) -- (11.250,-4.850);
\draw[wt/edge] (11.750,-4.850) -- (14.500,-4.850);
\draw[wt/bracket,draw=wtblue] (8.500,-4.100) -- (8.500,-3.980) -- (14.500,-3.980) -- (14.500,-4.100);
\draw[wt/bracket,draw=wtblue] (8.500,-5.730) -- (8.500,-5.850) -- (14.500,-5.850) -- (14.500,-5.730);
\draw[draw=wtblack,fill=white,line width=.70pt] (2.200,0.000) circle[radius=.073cm];
\draw[draw=wtblack,fill=white,line width=.70pt] (6.700,0.000) circle[radius=.073cm];
\draw[draw=wtblack,fill=white,line width=.70pt] (8.500,0.000) circle[radius=.073cm];
\draw[draw=wtblack,fill=white,line width=.70pt] (14.500,0.000) circle[radius=.073cm];
\draw[wt/edge,fill=white] (0.650,0.820) circle[radius=.073cm];
\draw[wt/edge,fill=white] (0.650,-0.820) circle[radius=.073cm];
\draw[draw=wtblack,fill=white,line width=.70pt] (2.200,-4.850) circle[radius=.073cm];
\draw[draw=wtblack,fill=white,line width=.70pt] (6.700,-4.850) circle[radius=.073cm];
\draw[draw=wtblack,fill=white,line width=.70pt] (8.500,-4.850) circle[radius=.073cm];
\draw[draw=wtblack,fill=white,line width=.70pt] (14.500,-4.850) circle[radius=.073cm];
\draw[wt/edge,fill=white] (0.650,-4.030) circle[radius=.073cm];
\draw[wt/edge,fill=white] (0.650,-5.670) circle[radius=.073cm];
\node[wt/math,text=wtblack,anchor=center] at (4.450,0.018) {$\displaystyle \cdots$};
\node[wt/math,text=wtblack,anchor=center] at (11.500,0.018) {$\displaystyle \cdots$};
\node[wt/math,text=wtblack,anchor=north west] at (2.280,-0.270) {$\displaystyle a_0=o$};
\node[wt/math,text=wtblack,anchor=north] at (6.700,-0.270) {$\displaystyle a_{k-1}$};
\node[wt/math,text=wtblack,anchor=north] at (8.500,-0.270) {$\displaystyle a_k$};
\node[wt/math,text=wtblack,anchor=north] at (14.500,-0.270) {$\displaystyle a_m$};
\node[wt/math,text=wtblack,anchor=east] at (0.330,0.820) {$\displaystyle b$};
\node[wt/math,text=wtblack,anchor=east] at (0.330,-0.820) {$\displaystyle c$};
\node[wt/math,text=wtblue,anchor=north] at (0.650,-1.140) {$\displaystyle Z$};
\node[wt/math,text=wtblue,anchor=center] at (8.600,1.610) {$\displaystyle F=\{b\}:\quad Z=\omega(c)$};
\node[wt/math,text=wtblue,anchor=south] at (4.600,0.550) {$\displaystyle S,\quad |S|=k+s$};
\node[wt/math,text=wtblue,anchor=south] at (11.500,1.060) {$\displaystyle T_{m+1-k}$};
\node[wt/math,text=wtblue,anchor=north] at (11.500,-1.190) {$\displaystyle Y_{m+1-k}$};
\node[wt/math,text=wtblue,anchor=center] at (8.100,-1.940) {$\displaystyle \sum_{v\notin S}\omega(v)=Y_{m+1-k}+Z$};
\node[wt/math,text=wtblack,anchor=center] at (4.450,-4.832) {$\displaystyle \cdots$};
\node[wt/math,text=wtblack,anchor=center] at (11.500,-4.832) {$\displaystyle \cdots$};
\node[wt/math,text=wtblack,anchor=north west] at (2.280,-5.120) {$\displaystyle a_0=o$};
\node[wt/math,text=wtblack,anchor=north] at (6.700,-5.120) {$\displaystyle a_{k-1}$};
\node[wt/math,text=wtblack,anchor=north] at (8.500,-5.120) {$\displaystyle a_k$};
\node[wt/math,text=wtblack,anchor=north] at (14.500,-5.120) {$\displaystyle a_m$};
\node[wt/math,text=wtblack,anchor=east] at (0.330,-4.030) {$\displaystyle b$};
\node[wt/math,text=wtblack,anchor=east] at (0.330,-5.670) {$\displaystyle c$};
\node[wt/math,text=wtblue,anchor=center] at (8.600,-3.240) {$\displaystyle F=\{b,c\}:\quad Z=0$};
\node[wt/math,text=wtblue,anchor=south] at (4.600,-4.300) {$\displaystyle S,\quad |S|=k+s$};
\node[wt/math,text=wtblue,anchor=south] at (11.500,-3.790) {$\displaystyle T_{m+1-k}$};
\node[wt/math,text=wtblue,anchor=north] at (11.500,-6.040) {$\displaystyle Y_{m+1-k}$};
\node[wt/math,text=wtblue,anchor=center] at (8.100,-6.790) {$\displaystyle \sum_{v\notin S}\omega(v)=Y_{m+1-k}+Z$};
\end{tikzpicture}
}
\caption{The path weight and the leaf weight $Z$ outside $S$.}
\label{fig:outside-mass}
\end{figure}

If $t\le m+1-k$, then $T_t\cap S=\varnothing$ and $Y_t'=Y_t$,
so \eqref{eq:two-leaf-tail-bound} follows from the induction hypothesis. If, in addition,
\[
Y_t'=\frac{t(2m+3)}{(m+2)^2},
\]
then $Y_t$ also attains equality before the operation. By \eqref{eq:two-leaf-tail-rigidity}, every
vertex of $T_t$ has weight $(2m+3)/(m+2)^2$ before the operation.
These weights do not change, so the same is true afterwards.

Now suppose that $t>m+1-k$, and define
\[
h=t-(m+1-k).
\]
Thus $h$ counts the vertices of $T_t\cap S$, and $1\le h\le k$.
We have
\[
T_t\setminus S=T_{m+1-k},
\qquad
T_t\cap S=\{a_{k-h},a_{k-h+1},\ldots,a_{k-1}\}.
\]
The operation leaves the weights on $T_{m+1-k}$ unchanged and assigns
weight $\mu$ to every vertex of $T_t\cap S$. Therefore,
\[
Y_t'=Y_{m+1-k}+h\mu.
\]
Substituting the expression for $\mu$ and rearranging gives
\begin{equation}\label{eq:two-leaf-deficit}
\begin{aligned}
\frac{t(2m+3)}{(m+2)^2}-Y_t'
={}&
\left(1-\frac{h}{k+s}\right)
\left(
\frac{(m+1-k)(2m+3)}{(m+2)^2}-Y_{m+1-k}
\right)\\
&+
\frac{h}{k+s}
\left(
\frac{(m+1+s)(2m+3)}{(m+2)^2}-2+Z
\right).
\end{aligned}
\end{equation}
Since $1\le h\le k<k+s$, the two coefficients of the bracketed terms
in \eqref{eq:two-leaf-deficit} are positive. The first bracketed term is nonnegative by \eqref{eq:two-leaf-tail-bound}.

If $s=1$, the second bracketed term is
\[
\frac{(m+2)(2m+3)}{(m+2)^2}-2+Z
=Z-\frac{1}{m+2}\ge0.
\]
If $s=2$, then $Z=0$ and the second bracketed term is
\[
\frac{(m+3)(2m+3)}{(m+2)^2}-2
=\frac{m+1}{(m+2)^2}>0.
\]
Thus the right-hand side of \eqref{eq:two-leaf-deficit} is nonnegative, proving \eqref{eq:two-leaf-tail-bound} for $\omega'$.

If
\[
Y_t'=\frac{t(2m+3)}{(m+2)^2},
\]
then both bracketed terms on the right-hand side of \eqref{eq:two-leaf-deficit} must vanish.
The preceding calculation shows that necessarily
\[
s=1,
\qquad
Z=\frac{1}{m+2},
\qquad
Y_{m+1-k}
=\frac{(m+1-k)(2m+3)}{(m+2)^2}.
\]
If $m+1-k>0$, the induction hypothesis \eqref{eq:two-leaf-tail-rigidity} implies that every vertex
of $T_{m+1-k}$ has weight $(2m+3)/(m+2)^2$, and the operation leaves
these weights unchanged. If $m+1-k=0$, then $T_{m+1-k}=\varnothing$.
Moreover,
\[
\mu=
\frac{
2-\dfrac{(m+1-k)(2m+3)}{(m+2)^2}
-\dfrac{1}{m+2}
}{k+1}
=\frac{2m+3}{(m+2)^2},
\]
so every vertex of $T_t\cap S$ also has weight $(2m+3)/(m+2)^2$ after
the operation. Since
\[
T_t=T_{m+1-k}\cup(T_t\cap S),
\]
we obtain \eqref{eq:two-leaf-tail-rigidity} for $\omega'$.

The induction is complete. Thus \eqref{eq:two-leaf-tail-bound}--\eqref{eq:two-leaf-tail-rigidity} hold after any finite sequence
of averaging operations. Taking $t=1$ in \eqref{eq:two-leaf-tail-bound} gives
\[
\omega(a_m)\le\frac{2m+3}{(m+2)^2}.
\]

We now determine when equality holds. Suppose that, after finitely many
averaging operations,
\[
\omega(a_m)=\frac{2m+3}{(m+2)^2}.
\]
Since the initial weight at $a_m$ is zero, some effective operation
contains $a_m$. Choose the last such operation, denote its selected set
by $S$, and write $\omega$ and $\omega'$ for the weights immediately before and
after it. No later effective operation contains $a_m$, so
\[
\omega'(a_m)=\frac{2m+3}{(m+2)^2}.
\]

If $S\cap\{b,c\}=\varnothing$, then connectivity and $a_m\in S$ imply
that, for some integer $0\le r\le m$,
\[
S=\{a_r,a_{r+1},\ldots,a_m\}.
\]
Put $t=m-r+1$, so that $S=T_t$. After the operation, every vertex of
$S$ has weight $(2m+3)/(m+2)^2$. Preservation of the total weight on
$S$ gives, before the operation,
\[
Y_t=\frac{t(2m+3)}{(m+2)^2}.
\]
By \eqref{eq:two-leaf-tail-rigidity}, the weights on $S$ are already equal before the operation,
contrary to its effectiveness. Therefore,
\[
S\cap\{b,c\}\ne\varnothing.
\]

Connectivity and $a_m\in S$ now imply
\[
\{a_0,a_1,\ldots,a_m\}\subseteq S.
\]
If both $b$ and $c$ belong to $S$, then $S=V(S(m,1,1))$, and hence
\[
\omega'(a_m)=\frac{2}{m+3}
<\frac{2m+3}{(m+2)^2}.
\]
The strict inequality follows from
\[
\frac{2m+3}{(m+2)^2}-\frac{2}{m+3}
=\frac{m+1}{(m+2)^2(m+3)}>0.
\]
This is a contradiction.

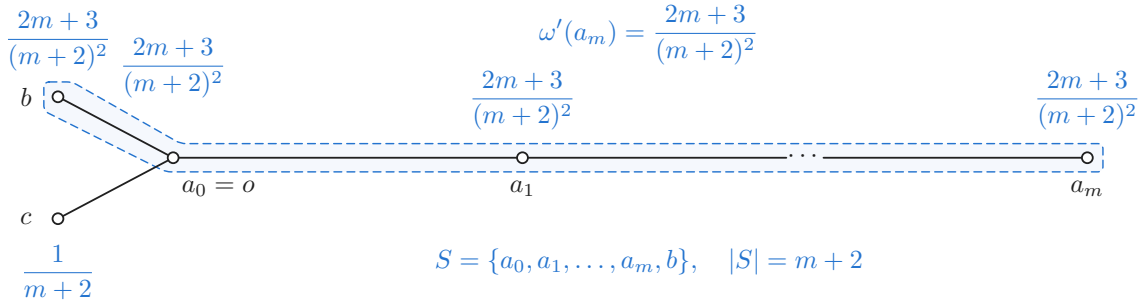
\begin{figure}[!htbp]
\centering
\resizebox{.96\linewidth}{!}{
\begin{tikzpicture}[wt/base]
\path (0,-2.0) rectangle (15.4,2.0);
\draw[wt/selection] (0.460,1.020) -- (0.840,1.020) -- (2.300,0.190) -- (14.710,0.190) -- (14.710,-0.190) -- (2.100,-0.190) -- (0.460,0.630) -- cycle;
\draw[wt/edge] (0.650,0.820) -- (2.200,0.000);
\draw[wt/edge] (0.650,-0.820) -- (2.200,0.000);
\draw[wt/edge] (2.200,0.000) -- (6.900,0.000);
\draw[wt/edge] (6.900,0.000) -- (10.450,0.000);
\draw[wt/edge] (10.950,0.000) -- (14.500,0.000);
\draw[draw=wtblack,fill=white,line width=.70pt] (2.200,0.000) circle[radius=.073cm];
\draw[draw=wtblack,fill=white,line width=.70pt] (6.900,0.000) circle[radius=.073cm];
\draw[draw=wtblack,fill=white,line width=.70pt] (14.500,0.000) circle[radius=.073cm];
\draw[wt/edge,fill=white] (0.650,0.820) circle[radius=.073cm];
\draw[wt/edge,fill=white] (0.650,-0.820) circle[radius=.073cm];
\node[wt/math,text=wtblack,anchor=center] at (10.700,0.018) {$\displaystyle \cdots$};
\node[wt/math,text=wtblue,anchor=south] at (2.200,0.780) {$\displaystyle \frac{2m+3}{(m+2)^2}$};
\node[wt/math,text=wtblack,anchor=north west] at (2.280,-0.270) {$\displaystyle a_0=o$};
\node[wt/math,text=wtblack,anchor=north] at (6.900,-0.270) {$\displaystyle a_1$};
\node[wt/math,text=wtblue,anchor=south] at (6.900,0.340) {$\displaystyle \frac{2m+3}{(m+2)^2}$};
\node[wt/math,text=wtblack,anchor=north] at (14.500,-0.270) {$\displaystyle a_m$};
\node[wt/math,text=wtblue,anchor=south] at (14.500,0.340) {$\displaystyle \frac{2m+3}{(m+2)^2}$};
\node[wt/math,text=wtblack,anchor=east] at (0.330,0.820) {$\displaystyle b$};
\node[wt/math,text=wtblue,anchor=south] at (0.650,1.140) {$\displaystyle \frac{2m+3}{(m+2)^2}$};
\node[wt/math,text=wtblack,anchor=east] at (0.330,-0.820) {$\displaystyle c$};
\node[wt/math,text=wtblue,anchor=north] at (0.650,-1.140) {$\displaystyle \frac1{m+2}$};
\node[wt/math,text=wtblue,anchor=center] at (8.600,1.650) {$\displaystyle \omega'(a_m)=\frac{2m+3}{(m+2)^2}$};
\node[wt/math,text=wtblue,anchor=center] at (8.600,-1.400) {$\displaystyle S=\{a_0,a_1,\ldots,a_m,b\},\quad |S|=m+2$};
\end{tikzpicture}
}
\caption{The equality configuration for two leaf sources, with leaf $b$ selected.}
\label{fig:double-leaf-equality}
\end{figure}

Thus $S$ contains exactly one of $b,c$. By the symmetry of the graph
and the initial weights in $b,c$, we may assume that
\[
S=\{a_0,a_1,\ldots,a_m,b\}.
\]
Then
\[
\omega'(a_m)=\frac{2-\omega(c)}{m+2}
=\frac{2m+3}{(m+2)^2},
\]
so
\[
\omega(c)=\frac{1}{m+2}.
\]
The operation does not select $c$, and therefore
\[
\omega'(c)=\frac{1}{m+2}.
\]
Every vertex of $S$ has weight $(2m+3)/(m+2)^2$ after the operation.
Interchanging $b$ and $c$ gives the other case.
One equality configuration is shown in Figure~\ref{fig:double-leaf-equality}.

Finally, in the case above, the initial weight at $c$ is one, whereas
its weight immediately before the last effective operation containing
$a_m$ is $1/(m+2)$. At least one earlier effective operation must
therefore have changed the weight at $c$. Together with the last
effective operation containing $a_m$, this gives at least two effective
operations in the sequence.
\end{proof}

\FloatBarrier
\begin{lemma}[Barrier for a single leaf source]
Suppose that the initial weight function $\omega$ on $S(m,1,1)$ satisfies
\[
\omega(a_i)=0 \quad (0\le i\le m),
\qquad
\omega(b)=1,
\qquad
\omega(c)=0.
\]
Then, after any finite sequence of averaging operations,
\[
\omega(a_m)\le\frac{1}{m+2}.
\]
If equality holds, then the last effective operation containing $a_m$
must average the set
\[
\{a_0,a_1,\ldots,a_m,b\}.
\]
Immediately after that operation,
\[
\omega(a_0)=\omega(a_1)=\cdots=\omega(a_m)=\omega(b)=\frac{1}{m+2},
\qquad
\omega(c)=0.
\]
Moreover, every sequence attaining equality contains at least one
effective operation.
\end{lemma}

\begin{proof}
We retain the notation $T_t$, where
\[
T_t=\{a_{m-t+1},a_{m-t+2},\ldots,a_m\}
\quad(1\le t\le m+1),
\qquad T_0=\varnothing.
\]
Thus $Y_t$ is the total weight on $T_t$.

Averaging preserves nonnegativity and the sum of all vertex weights.
Hence, throughout the sequence,
\[
\omega(v)\ge0
\quad\bigl(v\in V(S(m,1,1))\bigr),
\qquad
\sum_{v\in V(S(m,1,1))}\omega(v)=1.
\]

We prove by induction on the number of operations that
\begin{equation}\label{eq:single-leaf-tail-bound}
Y_t\le\frac{t}{m+2},
\qquad 0\le t\le m+1,
\end{equation}
and, simultaneously, that for $1\le t\le m+1$,
\begin{equation}\label{eq:single-leaf-tail-rigidity}
Y_t=\frac{t}{m+2}
\quad\Longrightarrow\quad
\omega(a_{m-t+1})=\cdots=\omega(a_m)=\frac{1}{m+2}.
\end{equation}
Taking $t=1$ in \eqref{eq:single-leaf-tail-bound} gives the required bound
\[
\omega(a_m)=Y_1\le\frac{1}{m+2}.
\]
We use \eqref{eq:single-leaf-tail-rigidity} to identify the last effective operation containing $a_m$
when equality holds.

Initially, every path vertex has weight zero, so
\[
Y_t=0,\qquad 0\le t\le m+1,
\]
and \eqref{eq:single-leaf-tail-bound} holds. For $1\le t\le m+1$, the equality in the antecedent of
\eqref{eq:single-leaf-tail-rigidity} does not hold initially, so \eqref{eq:single-leaf-tail-rigidity} holds as well.

Suppose that the current weight function $\omega$ satisfies \eqref{eq:single-leaf-tail-bound} and \eqref{eq:single-leaf-tail-rigidity}.
Let $S$ be the connected set selected by the next operation, let $\omega'$
be the resulting weight function, and put
\[
Y_t'=\sum_{v\in T_t}\omega'(v),
\qquad 0\le t\le m+1.
\]
Thus $Y_t'$ is the total weight on $T_t$ after the operation, with
$Y_0'=0$.

If the operation is ineffective, then $\omega'=\omega$ and the conclusion follows.
We may therefore assume that it is effective.

First suppose that
\[
S\cap\{b,c\}=\varnothing.
\]
Since $S$ is connected, there are integers $0\le r\le j\le m$ such that
\[
S=\{a_r,a_{r+1},\ldots,a_j\}.
\]
Apply the interpolation identity \eqref{eq:path-interpolation} and its equality analysis from the
proof of Lemma~1, using the induction hypotheses \eqref{eq:single-leaf-tail-bound} and \eqref{eq:single-leaf-tail-rigidity}.
It follows that $\omega'$ still satisfies \eqref{eq:single-leaf-tail-bound} and \eqref{eq:single-leaf-tail-rigidity}.

Next suppose that
\[
S\cap\{b,c\}\ne\varnothing.
\]
Define
\[
F=S\cap\{b,c\},
\qquad
s=|F|,
\qquad
k=\bigl|S\cap\{a_0,a_1,\ldots,a_m\}\bigr|.
\]
Thus $F$ is the set of selected vertices among $b,c$, $s$ is their
number, and $k$ is the number of selected vertices on the path
$\langle a_0,a_1,\ldots,a_m\rangle$.
Since the operation is effective and $S$ is connected,
\[
s\in\{1,2\},
\qquad
1\le k\le m+1,
\qquad
S=\{a_0,a_1,\ldots,a_{k-1}\}\cup F.
\]

Define
\[
Z=\sum_{u\in\{b,c\}\setminus F}\omega(u).
\]
Thus $Z$ is the total weight on the unselected vertices among $b,c$.
In particular, $Z\ge0$, and $Z=0$ when $s=2$.
Let $\mu$ be the common weight on $S$ after the operation.
Since the vertices outside $S$ consist precisely of $T_{m+1-k}$ and
$\{b,c\}\setminus F$,
\[
\mu=\frac{1-Y_{m+1-k}-Z}{k+s}.
\]

If $t\le m+1-k$, then $T_t\cap S=\varnothing$, so
\[
Y_t'=Y_t\le\frac{t}{m+2}.
\]
If $Y_t'=t/(m+2)$, then $Y_t$ also attains equality before the operation.
By \eqref{eq:single-leaf-tail-rigidity}, every vertex of $T_t$ has weight $1/(m+2)$ before the operation.
These weights do not change, so the same is true afterwards.

Now suppose that $t>m+1-k$, and define
\[
h=t-(m+1-k).
\]
Thus $h$ counts the vertices of $T_t\cap S$, and $1\le h\le k$.
We have
\[
T_t\setminus S=T_{m+1-k},
\qquad
T_t\cap S=\{a_{k-h},a_{k-h+1},\ldots,a_{k-1}\}.
\]
The operation leaves all weights on $T_{m+1-k}$ unchanged and assigns
weight $\mu$ to each vertex of $T_t\cap S$. Therefore,
\[
Y_t'=Y_{m+1-k}+h\mu.
\]
Substituting the expression for $\mu$ and rearranging gives
\begin{equation}\label{eq:single-leaf-deficit}
\begin{aligned}
\frac{t}{m+2}-Y_t'
={}&
\left(1-\frac{h}{k+s}\right)
\left(\frac{m+1-k}{m+2}-Y_{m+1-k}\right)\\
&+
\frac{h}{k+s}
\left(\frac{m+1+s}{m+2}-1+Z\right).
\end{aligned}
\end{equation}
Since $1\le h\le k<k+s$, the two coefficients of the bracketed terms
are positive. The first bracketed term is nonnegative by \eqref{eq:single-leaf-tail-bound}.
If $s=1$, the second bracketed term is $Z\ge0$.
If $s=2$, then $Z=0$ and the second bracketed term is
\[
\frac{m+3}{m+2}-1=\frac{1}{m+2}>0.
\]
Thus the right-hand side of \eqref{eq:single-leaf-deficit} is nonnegative, proving \eqref{eq:single-leaf-tail-bound} for $\omega'$.

If $Y_t'=t/(m+2)$, then both bracketed terms on the right-hand side of
\eqref{eq:single-leaf-deficit} must vanish. Hence
\[
s=1,
\qquad
Z=0,
\qquad
Y_{m+1-k}=\frac{m+1-k}{m+2}.
\]
If $m+1-k>0$, the induction hypothesis \eqref{eq:single-leaf-tail-rigidity} implies that every vertex of
$T_{m+1-k}$ has weight $1/(m+2)$, and the operation leaves these weights
unchanged. If $m+1-k=0$, then $T_{m+1-k}=\varnothing$. Moreover,
\[
\mu=
\frac{1-\dfrac{m+1-k}{m+2}}{k+1}
=\frac{1}{m+2},
\]
so every vertex of $T_t\cap S$ also has weight $1/(m+2)$ after the operation.
Since
\[
T_t=T_{m+1-k}\cup(T_t\cap S),
\]
we obtain \eqref{eq:single-leaf-tail-rigidity} for $\omega'$.

The induction is complete. Thus \eqref{eq:single-leaf-tail-bound} and \eqref{eq:single-leaf-tail-rigidity} hold after any finite
sequence of averaging operations. Taking $t=1$ in \eqref{eq:single-leaf-tail-bound} gives
\[
\omega(a_m)\le\frac{1}{m+2}.
\]

We now determine when equality holds. Suppose that, after finitely many
averaging operations,
\[
\omega(a_m)=\frac{1}{m+2}.
\]
Since the initial weight at $a_m$ is zero, some effective operation
contains $a_m$. Choose the last such operation, denote its selected set
by $S$, and write $\omega$ and $\omega'$ for the weights immediately before and
after it. No later effective operation contains $a_m$, so
\[
\omega'(a_m)=\frac{1}{m+2}.
\]

If $S\cap\{b,c\}=\varnothing$, then connectivity and $a_m\in S$ imply
that, for some integer $0\le r\le m$,
\[
S=\{a_r,a_{r+1},\ldots,a_m\}.
\]
Put $t=m-r+1$, so that $S=T_t$. After the operation, every vertex of
$S$ has weight $1/(m+2)$. Preservation of the total weight on $S$ gives,
before the operation,
\[
Y_t=\frac{t}{m+2}.
\]
By \eqref{eq:single-leaf-tail-rigidity}, the weights on $S$ are already equal before the operation,
contrary to its effectiveness. Therefore,
\[
S\cap\{b,c\}\ne\varnothing.
\]
Connectivity and $a_m\in S$ further imply
\[
\{a_0,a_1,\ldots,a_m\}\subseteq S.
\]

If both $b$ and $c$ belong to $S$, then $S=V(S(m,1,1))$, so
\[
\omega'(a_m)=\frac{1}{m+3}<\frac{1}{m+2}.
\]
This is a contradiction. Thus $S$ contains exactly one of $b,c$.

Observe that the weight at $b$ remains positive after any finite sequence
of operations. Indeed, its initial weight is one. An operation avoiding
$b$ leaves its weight unchanged. If an operation contains $b$ and
$\omega(b)>0$ before the operation, then nonnegativity makes the total weight
on the selected set positive, so the new weight at $b$ is again positive.
The assertion follows by induction on the number of operations.

\begin{figure}[!htbp]
\centering
\resizebox{.96\linewidth}{!}{
\begin{tikzpicture}[wt/base]
\path (0,-2.0) rectangle (15.4,2.0);
\draw[wt/selection] (0.460,1.020) -- (0.840,1.020) -- (2.300,0.190) -- (14.710,0.190) -- (14.710,-0.190) -- (2.100,-0.190) -- (0.460,0.630) -- cycle;
\draw[wt/edge] (0.650,0.820) -- (2.200,0.000);
\draw[wt/edge] (0.650,-0.820) -- (2.200,0.000);
\draw[wt/edge] (2.200,0.000) -- (6.900,0.000);
\draw[wt/edge] (6.900,0.000) -- (10.450,0.000);
\draw[wt/edge] (10.950,0.000) -- (14.500,0.000);
\draw[draw=wtblack,fill=white,line width=.70pt] (2.200,0.000) circle[radius=.073cm];
\draw[draw=wtblack,fill=white,line width=.70pt] (6.900,0.000) circle[radius=.073cm];
\draw[draw=wtblack,fill=white,line width=.70pt] (14.500,0.000) circle[radius=.073cm];
\draw[wt/edge,fill=white] (0.650,0.820) circle[radius=.073cm];
\draw[wt/edge,fill=white] (0.650,-0.820) circle[radius=.073cm];
\node[wt/math,text=wtblack,anchor=center] at (10.700,0.018) {$\displaystyle \cdots$};
\node[wt/math,text=wtblue,anchor=south] at (2.200,0.780) {$\displaystyle \frac1{m+2}$};
\node[wt/math,text=wtblack,anchor=north west] at (2.280,-0.270) {$\displaystyle a_0=o$};
\node[wt/math,text=wtblack,anchor=north] at (6.900,-0.270) {$\displaystyle a_1$};
\node[wt/math,text=wtblue,anchor=south] at (6.900,0.340) {$\displaystyle \frac1{m+2}$};
\node[wt/math,text=wtblack,anchor=north] at (14.500,-0.270) {$\displaystyle a_m$};
\node[wt/math,text=wtblue,anchor=south] at (14.500,0.340) {$\displaystyle \frac1{m+2}$};
\node[wt/math,text=wtblack,anchor=east] at (0.330,0.820) {$\displaystyle b$};
\node[wt/math,text=wtblue,anchor=south] at (0.650,1.140) {$\displaystyle \frac1{m+2}$};
\node[wt/math,text=wtblack,anchor=east] at (0.330,-0.820) {$\displaystyle c$};
\node[wt/math,text=wtblue,anchor=north] at (0.650,-1.140) {$\displaystyle 0$};
\node[wt/math,text=wtblue,anchor=center] at (8.600,1.650) {$\displaystyle \omega'(a_m)=\frac1{m+2}$};
\node[wt/math,text=wtblue,anchor=center] at (8.600,-1.400) {$\displaystyle S=\{a_0,a_1,\ldots,a_m,b\},\quad |S|=m+2$};
\end{tikzpicture}
}
\caption{The equality configuration for a single leaf source.}
\label{fig:single-leaf-equality}
\end{figure}
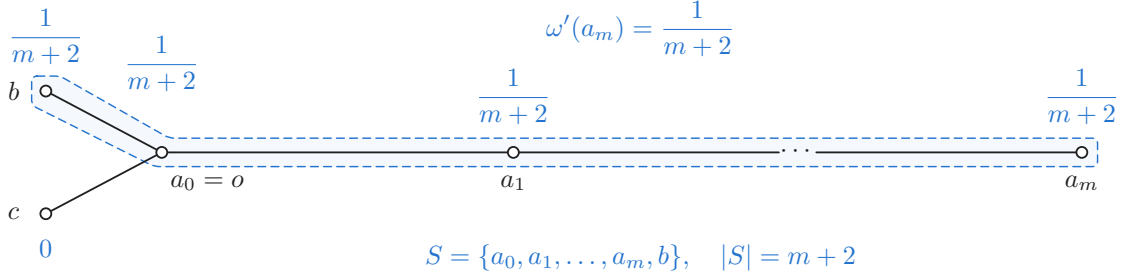

Consequently, if
\[
S=\{a_0,a_1,\ldots,a_m,c\},
\]
then $\omega(b)>0$ before the operation, and hence
\[
\omega'(a_m)=\frac{1-\omega(b)}{m+2}<\frac{1}{m+2},
\]
again a contradiction. We must therefore have
\[
S=\{a_0,a_1,\ldots,a_m,b\}.
\]
In this case,
\[
\omega'(a_m)=\frac{1-\omega(c)}{m+2}=\frac{1}{m+2},
\]
so $\omega(c)=0$. The operation does not select $c$, so $\omega'(c)=0$, while
every vertex of $S$ has weight $1/(m+2)$ after the operation. Thus,
immediately after the operation,
\[
\omega'(a_0)=\omega'(a_1)=\cdots=\omega'(a_m)=\omega'(b)=\frac{1}{m+2},
\qquad
\omega'(c)=0.
\]
This equality configuration is shown in Figure~\ref{fig:single-leaf-equality}.

Finally, the weight at $a_m$ changes from zero to a positive value.
Every sequence attaining equality therefore contains at least one
effective operation.
\end{proof}

\FloatBarrier
\section{Construction of the reference profile}
\label{sec:reference-profile}
We first express the target weight produced by a forward strategy as the
inner product of a vector with the initial weights. We then construct a
reference vector using $m+2$ averaging operations and compute its
coordinates. This construction will allow us to choose nonnegative initial
weights for which applying these operations in reverse order produces a
strictly larger target weight than any strategy of at most $m+1$
operations. It will follow that an optimal strategy requires at least
$m+2$ operations.

Fix an integer $m\ge2$, put $G=S(m,1,1)$, and retain the convention
$a_0=o$. We identify a vertex weight function with a vector indexed by
$V(G)$; its coordinate at $v$ is the weight of $v$.
To introduce the reverse profiles used below, we first describe the
relation between averaging operators and the target weight.

For a nonempty connected set $S\subseteq V(G)$, let $A_S$ be the
averaging operator on $S$. For a real-valued weight function
$w:V(G)\to\mathbb R$, define
\[
(A_Sw)(v)=
\begin{cases}
\displaystyle\frac1{|S|}\sum_{x\in S}w(x), & v\in S,\\[6pt]
w(v), & v\notin S.
\end{cases}
\]
We allow real weights here to accommodate the transformations used later.
For real-valued weight functions $f,g$, define the inner product by
$\langle f,g\rangle=\sum_{v\in V(G)}f(v)g(v)$.
A direct calculation gives
\[
\begin{aligned}
\langle A_Sf,g\rangle
&=\sum_{v\notin S}f(v)g(v)
+\frac1{|S|}
 \left(\sum_{v\in S}f(v)\right)
 \left(\sum_{v\in S}g(v)\right)\\
&=\langle f,A_Sg\rangle.
\end{aligned}
\]

\begin{definition}[Reverse profile]
Let $e_o$ be the unit vector at $o$, that is,
\[
e_o(v)=
\begin{cases}
1, & v=o,\\
0, & v\ne o.
\end{cases}
\]
Suppose that a forward strategy averages the nonempty connected sets
$S_1,\ldots,S_r$ in this order, where $r$ is a nonnegative integer.
The vector obtained from $e_o$ by averaging these sets in reverse order,
$S_r,\ldots,S_1$, namely
\[
q=A_{S_1}\cdots A_{S_r}e_o
\]
is called the \emph{reverse profile} of the strategy.
Operator products act from right to left. When $r=0$, no operation is
performed and we set $q=e_o$.
\end{definition}

By the inner product identity above, the target weight at $o$ obtained by
this forward strategy from initial weights $\omega$ is
\begin{equation}\label{eq:forward-reverse}
\begin{aligned}
\bigl(A_{S_r}\cdots A_{S_1}\omega\bigr)(o)
&=\langle e_o,A_{S_r}\cdots A_{S_1}\omega\rangle\\
&=\langle A_{S_1}\cdots A_{S_r}e_o,\omega\rangle\\
&=\langle q,\omega\rangle
=\sum_{v\in V(G)}q(v)\omega(v).
\end{aligned}
\end{equation}
Thus $q(v)$ is the coefficient of the initial weight $\omega(v)$ in the
final target weight. The vector $q$ depends only on the graph, the target
vertex and the sequence of operations, and is independent of $\omega$.

\begin{definition}[Reachable profiles and minimum length]
For a nonnegative integer $k$, let $\mathcal R_k$ be the set of weight
vectors reachable from $e_o$ by at most $k$ averaging operations.
Each such vector is the reverse profile of a forward strategy.
In particular, $\mathcal R_0=\{e_o\}$, while $\mathcal R_{m+1}$ and
$\mathcal R_{m+2}$ correspond to at most $m+1$ and at most $m+2$
operations, respectively. Put
$\mathcal R_{\mathrm{fin}}=\bigcup_{k\ge0}\mathcal R_k$, the set of
all profiles reachable from $e_o$ by finitely many averaging operations.
For $q\in\mathcal R_{\mathrm{fin}}$, define
$\ell(q)=\min\{k\ge0:q\in\mathcal R_k\}$, where $k$ ranges over the
nonnegative integers. Thus $\ell(q)$ is the minimum number of averaging
operations needed to obtain $q$ from $e_o$.
\end{definition}

Averaging preserves nonnegativity and the sum of all coordinates.
Since the coordinates of $e_o$ sum to one, every
$q\in\mathcal R_{\mathrm{fin}}$ satisfies
\[
0\le q(v)\le1\quad(v\in V(G)),
\qquad
\sum_{v\in V(G)}q(v)=1.
\]
An ineffective operation leaves the vector unchanged. Hence every sequence
of length $\ell(q)$ producing $q$ consists entirely of effective operations.

By \eqref{eq:forward-reverse}, to rule out all optimal strategies of at
most $m+1$ operations, it suffices to construct
$p_m\in\mathcal R_{m+2}$ and nonnegative initial weights $\omega_m$ such that
\[
\langle q,\omega_m\rangle<\langle p_m,\omega_m\rangle
\qquad(q\in\mathcal R_{m+1}).
\]
We construct $p_m$ in this section. The next two sections compare its
coordinates with other reachable profiles and choose $\omega_m$.

For $0\le j\le m$, put $P_j=\{a_0,a_1,\ldots,a_j\}$, the first
$j+1$ vertices of $\langle a_0,a_1,\ldots,a_m\rangle$ starting at $a_0$.
Starting from $e_o$, first average the sets
\[
P_m,\qquad P_{m-1}\cup\{b\},\qquad P_{m-1}\cup\{c\}
\]
in this order. Then, for $j=m-2,m-3,\ldots,0$, successively average
\[
\begin{cases}
P_j\cup\{b\}, & m-j\text{ is even},\\
P_j\cup\{c\}, & m-j\text{ is odd}.
\end{cases}
\]
All these sets are connected. There are $1+2+(m-1)=m+2$ operations in
total. Denote the resulting vector by $p_m$; thus
$p_m\in\mathcal R_{m+2}$. After the third operation, each successive
selected set contains one fewer path vertex and alternates between the
leaves $b$ and $c$.

\begin{figure}[!htbp]
\centering
\resizebox{.96\linewidth}{!}{
\begin{tikzpicture}[wt/base]
\path (0,-9.8) rectangle (15.4,1.8);
\draw[wt/selection] (1.990,-0.190) rectangle (14.710,0.190);
\draw[wt/selection] (0.460,-2.880) -- (0.840,-2.880) -- (2.300,-3.710) -- (9.110,-3.710) -- (9.110,-4.090) -- (2.100,-4.090) -- (0.460,-3.270) -- cycle;
\draw[wt/selection] (0.460,-8.820) -- (0.840,-8.820) -- (2.300,-7.990) -- (9.110,-7.990) -- (9.110,-7.610) -- (2.100,-7.610) -- (0.460,-8.430) -- cycle;
\draw[wt/edge] (0.650,0.820) -- (2.200,0.000);
\draw[wt/edge] (0.650,-0.820) -- (2.200,0.000);
\draw[wt/edge] (2.200,0.000) -- (5.300,0.000);
\draw[wt/edge] (5.800,0.000) -- (8.900,0.000);
\draw[wt/edge] (8.900,0.000) -- (14.500,0.000);
\draw[wt/edge] (0.650,-3.080) -- (2.200,-3.900);
\draw[wt/edge] (0.650,-4.720) -- (2.200,-3.900);
\draw[wt/edge] (2.200,-3.900) -- (5.300,-3.900);
\draw[wt/edge] (5.800,-3.900) -- (8.900,-3.900);
\draw[wt/edge] (8.900,-3.900) -- (14.500,-3.900);
\draw[wt/edge] (0.650,-6.980) -- (2.200,-7.800);
\draw[wt/edge] (0.650,-8.620) -- (2.200,-7.800);
\draw[wt/edge] (2.200,-7.800) -- (5.300,-7.800);
\draw[wt/edge] (5.800,-7.800) -- (8.900,-7.800);
\draw[wt/edge] (8.900,-7.800) -- (14.500,-7.800);
\draw[draw=wtblack,fill=white,line width=.70pt] (2.200,0.000) circle[radius=.073cm];
\draw[draw=wtblack,fill=white,line width=.70pt] (8.900,0.000) circle[radius=.073cm];
\draw[draw=wtblack,fill=white,line width=.70pt] (14.500,0.000) circle[radius=.073cm];
\draw[wt/edge,fill=white] (0.650,0.820) circle[radius=.073cm];
\draw[wt/edge,fill=white] (0.650,-0.820) circle[radius=.073cm];
\draw[draw=wtblack,fill=white,line width=.70pt] (2.200,-3.900) circle[radius=.073cm];
\draw[draw=wtblack,fill=white,line width=.70pt] (8.900,-3.900) circle[radius=.073cm];
\draw[draw=wtblack,fill=white,line width=.70pt] (14.500,-3.900) circle[radius=.073cm];
\draw[wt/edge,fill=white] (0.650,-3.080) circle[radius=.073cm];
\draw[wt/edge,fill=white] (0.650,-4.720) circle[radius=.073cm];
\draw[draw=wtblack,fill=white,line width=.70pt] (2.200,-7.800) circle[radius=.073cm];
\draw[draw=wtblack,fill=white,line width=.70pt] (8.900,-7.800) circle[radius=.073cm];
\draw[draw=wtblack,fill=white,line width=.70pt] (14.500,-7.800) circle[radius=.073cm];
\draw[wt/edge,fill=white] (0.650,-6.980) circle[radius=.073cm];
\draw[wt/edge,fill=white] (0.650,-8.620) circle[radius=.073cm];
\node[wt/math,text=wtblack,anchor=center] at (5.550,0.018) {$\displaystyle \cdots$};
\node[wt/math,text=wtblue,anchor=south] at (2.200,0.340) {$\displaystyle \frac1{m+1}$};
\node[wt/math,text=wtblack,anchor=north west] at (2.280,-0.270) {$\displaystyle a_0=o$};
\node[wt/math,text=wtblack,anchor=north] at (8.900,-0.270) {$\displaystyle a_{m-1}$};
\node[wt/math,text=wtblue,anchor=south] at (8.900,0.340) {$\displaystyle \frac1{m+1}$};
\node[wt/math,text=wtblack,anchor=north] at (14.500,-0.270) {$\displaystyle a_m$};
\node[wt/math,text=wtblue,anchor=south] at (14.500,0.340) {$\displaystyle \frac1{m+1}$};
\node[wt/math,text=wtblack,anchor=east] at (0.330,0.820) {$\displaystyle b$};
\node[wt/math,text=wtblue,anchor=south] at (0.650,1.140) {$\displaystyle 0$};
\node[wt/math,text=wtblack,anchor=east] at (0.330,-0.820) {$\displaystyle c$};
\node[wt/math,text=wtblue,anchor=north] at (0.650,-1.140) {$\displaystyle 0$};
\node[wt/math,text=wtblue,anchor=center] at (9.000,1.550) {$\displaystyle 1:\quad P_m$};
\node[wt/math,text=wtblack,anchor=center] at (5.550,-3.882) {$\displaystyle \cdots$};
\node[wt/math,text=wtblue,anchor=south] at (2.200,-3.340) {$\displaystyle \rho_m$};
\node[wt/math,text=wtblack,anchor=north west] at (2.280,-4.170) {$\displaystyle a_0=o$};
\node[wt/math,text=wtblack,anchor=north] at (8.900,-4.170) {$\displaystyle a_{m-1}$};
\node[wt/math,text=wtblue,anchor=south] at (8.900,-3.560) {$\displaystyle \rho_m$};
\node[wt/math,text=wtblack,anchor=north] at (14.500,-4.170) {$\displaystyle a_m$};
\node[wt/math,text=wtblue,anchor=south] at (14.500,-3.560) {$\displaystyle \frac1{m+1}$};
\node[wt/math,text=wtblack,anchor=east] at (0.330,-3.080) {$\displaystyle b$};
\node[wt/math,text=wtblue,anchor=south] at (0.650,-2.760) {$\displaystyle \rho_m$};
\node[wt/math,text=wtblack,anchor=east] at (0.330,-4.720) {$\displaystyle c$};
\node[wt/math,text=wtblue,anchor=north] at (0.650,-5.040) {$\displaystyle 0$};
\node[wt/math,text=wtblue,anchor=center] at (9.000,-2.350) {$\displaystyle 2:\quad P_{m-1}\cup\{b\}$};
\node[wt/math,text=wtblack,anchor=center] at (5.550,-7.782) {$\displaystyle \cdots$};
\node[wt/math,text=wtblue,anchor=south] at (2.200,-7.460) {$\displaystyle \rho_{m-1}$};
\node[wt/math,text=wtblack,anchor=north west] at (2.280,-8.070) {$\displaystyle a_0=o$};
\node[wt/math,text=wtblack,anchor=north] at (8.900,-8.070) {$\displaystyle a_{m-1}$};
\node[wt/math,text=wtblue,anchor=south] at (8.900,-7.460) {$\displaystyle \rho_{m-1}$};
\node[wt/math,text=wtblack,anchor=north] at (14.500,-8.070) {$\displaystyle a_m$};
\node[wt/math,text=wtblue,anchor=south] at (14.500,-7.460) {$\displaystyle \frac1{m+1}$};
\node[wt/math,text=wtblack,anchor=east] at (0.330,-6.980) {$\displaystyle b$};
\node[wt/math,text=wtblue,anchor=south] at (0.650,-6.660) {$\displaystyle \rho_m$};
\node[wt/math,text=wtblack,anchor=east] at (0.330,-8.620) {$\displaystyle c$};
\node[wt/math,text=wtblue,anchor=north] at (0.650,-8.940) {$\displaystyle \rho_{m-1}$};
\node[wt/math,text=wtblue,anchor=center] at (9.000,-6.250) {$\displaystyle 3:\quad P_{m-1}\cup\{c\}$};
\end{tikzpicture}
}
\caption{The first three operations constructing $p_m$ and their resulting weights.}
\label{fig:reference-start}
\end{figure}
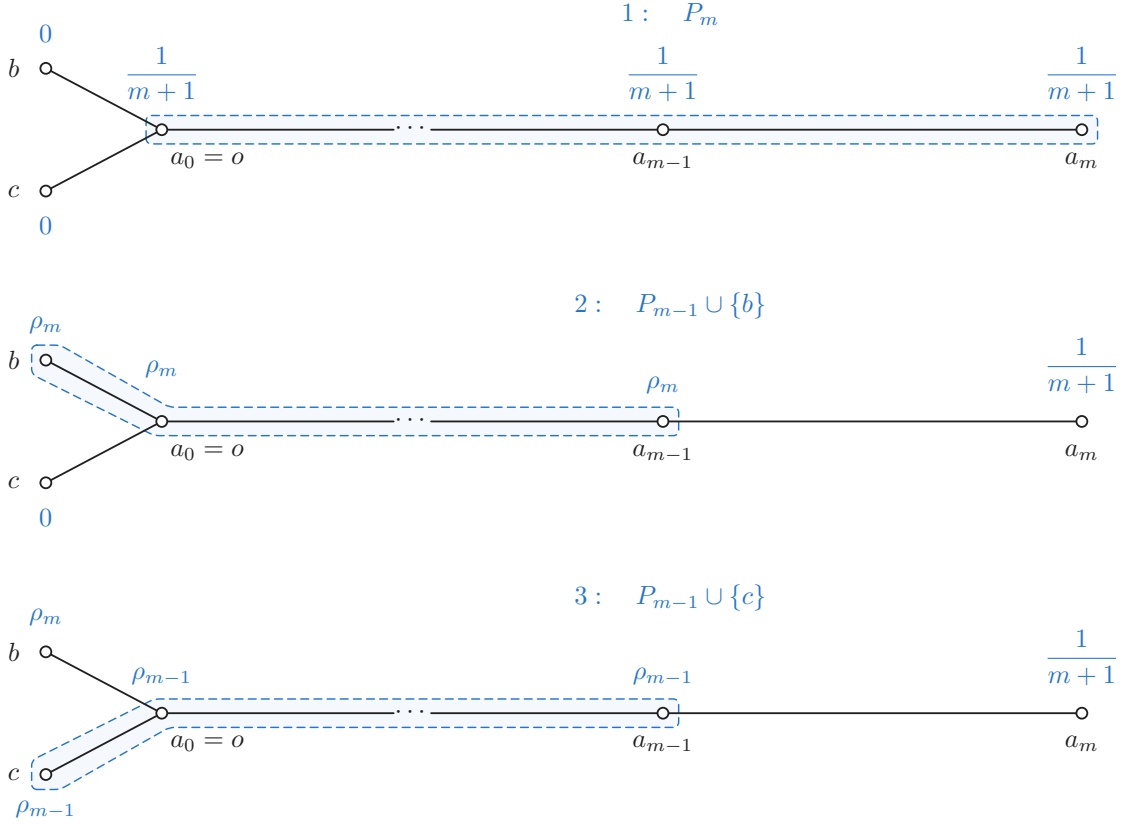

These operations successively fix the final weights at
$a_m,a_{m-1},\ldots,a_0$. After the first operation, every vertex of
$P_m$ has weight $1/(m+1)$, and no later operation contains $a_m$.
Write the common weights on the selected sets after the second and
third operations as
\[
\rho_m=\frac{m}{(m+1)^2},
\qquad
\rho_{m-1}=\frac{m^2}{(m+1)^3},
\]
respectively. Here $\rho_m$ is the intermediate weight produced by the
second operation, rather than $p_m(a_m)$. After the third operation,
every vertex of $P_{m-1}\cup\{c\}$ has weight $\rho_{m-1}$, while $b$
has weight $\rho_m$. The first three selected sets and the resulting
weights are shown in Figure~\ref{fig:reference-start}.

At a later stage, suppose that every vertex of $P_j$ and one of $b,c$
has weight $\rho_j$, while the other leaf has weight $\rho_{j+1}$.
The next operation averages $P_{j-1}$ together with the leaf of weight
$\rho_{j+1}$. Write the resulting common weight as $\rho_{j-1}$.
Then
\begin{equation}\label{eq:reference-recursion}
\rho_{j-1}=\frac{j\rho_j+\rho_{j+1}}{j+1},
\qquad 1\le j\le m-1.
\end{equation}
The vertex $a_j$ is excluded from this and all later operations, so its
final weight is $\rho_j$. This recursively determines
$\rho_{m-2},\ldots,\rho_0$ and gives
\[
p_m(a_m)=\frac1{m+1},
\qquad
p_m(a_j)=\rho_j\quad(0\le j\le m-1),
\]
and
\[
\bigl(p_m(b),p_m(c)\bigr)=
\begin{cases}
(\rho_0,\rho_1), & m\text{ is even},\\
(\rho_1,\rho_0), & m\text{ is odd}.
\end{cases}
\]

The recurrence also determines the relative sizes of consecutive
$\rho_j$. First, $\rho_m-\rho_{m-1}=m/(m+1)^3>0$.
By \eqref{eq:reference-recursion},
\[
\rho_j-\rho_{j-1}
=-\frac{\rho_{j+1}-\rho_j}{j+1},
\qquad 1\le j\le m-1.
\]
Thus all differences $\rho_{j+1}-\rho_j$ are nonzero, with alternating signs:
\begin{equation}\label{eq:reference-signs}
(-1)^{m-1-j}(\rho_{j+1}-\rho_j)>0,
\qquad 0\le j\le m-1.
\end{equation}

\FloatBarrier
\section{Successive coordinate comparisons of reachable profiles}

We compare $q\in\mathcal R_{\mathrm{fin}}$ with $p_m$ in the order
$a_m,a_{m-1},\ldots,a_0$. We proceed to the next coordinate only if all
previous coordinates are equal; the comparison stops at the first
unequal coordinate. The three barrier lemmas in Section~2 provide the
coordinate inequalities. Their equality conditions determine the
resulting configurations and the number of effective operations required.

During the comparison, we transform the weights to turn the equality
configuration from the previous step into one of the initial configurations
in Section~2: all remaining path vertices have weight zero, and the
weights at $b,c$ are either both one or one and zero. We can then apply
the barrier for two leaf sources or for a single leaf source and
translate the resulting inequality back to the original weights.
For any real $\alpha$, nonzero real $\beta$, and real-valued weight
function $w$,
\begin{equation}\label{eq:affine-average}
\frac{(A_Sw)(v)-\alpha}{\beta}=
\begin{cases}
\displaystyle\frac1{|S|}
\sum_{x\in S}\frac{w(x)-\alpha}{\beta}, & v\in S,\\[6pt]
\displaystyle\frac{w(v)-\alpha}{\beta}, & v\notin S.
\end{cases}
\end{equation}

Thus applying $w(v)\mapsto(w(v)-\alpha)/\beta$ to every coordinate and
then averaging gives the same result as averaging first and applying the
transformation afterwards. Since $\beta\ne0$, two coordinates are equal
before the transformation if and only if they are equal afterwards.

To express the comparisons uniformly, define
\begin{equation}\label{eq:alternating-profile}
\Lambda_m(q)=
\bigl(q(a_m),-q(a_{m-1}),q(a_{m-2}),\ldots,(-1)^mq(a_0)\bigr).
\end{equation}

The coordinate at position $r+1$ is $(-1)^rq(a_{m-r})$, where $0\le r\le m$.
We use the lexicographic order $\le_{\mathrm{lex}}$: the first unequal
coordinate, reading from left to right, determines the order, and the
vectors are equal if all coordinates agree. For example, if
$q(a_m)=p_m(a_m)$ and $q(a_{m-1})>p_m(a_{m-1})$, then
$-q(a_{m-1})<-p_m(a_{m-1})$, giving
$\Lambda_m(q)<_{\mathrm{lex}}\Lambda_m(p_m)$.
The minus signs in \eqref{eq:alternating-profile} turn comparisons of the
form $q(a_i)\ge p_m(a_i)$ into $-q(a_i)\le-p_m(a_i)$.
The direction at each step is determined by \eqref{eq:reference-signs},
as computed in the proof below.

\begin{theorem}[Lexicographic comparison and a lower bound on length]
\label{thm:profile-comparison}
For every $q\in\mathcal R_{\mathrm{fin}}$,
\[
\Lambda_m(q)\le_{\mathrm{lex}}\Lambda_m(p_m).
\]
If $\Lambda_m(q)=\Lambda_m(p_m)$, then $\ell(q)\ge m+2$.
For every $q\in\mathcal R_{m+1}$,
\[
\Lambda_m(q)<_{\mathrm{lex}}\Lambda_m(p_m).
\]
\end{theorem}

\begin{proof}
Fix $q\in\mathcal R_{\mathrm{fin}}$ and put $k=\ell(q)$.
Choose a shortest sequence producing $q$ from $e_o$, and denote its
successively averaged nonempty connected sets by $S_1,\ldots,S_k$.
Every operation in this sequence is effective. Put
\[
q^{(0)}=e_o,\qquad
q^{(s)}=A_{S_s}q^{(s-1)}\quad(1\le s\le k).
\]
Then $q^{(k)}=q$. The superscript $s$ indicates the time immediately
after the $s$th operation, and $q^{(s)}(v)$ is the weight at $v$ at that time.

We first compare the coordinate at $a_m$. By the root-source barrier,
\[
q(a_m)\le\frac1{m+1}=p_m(a_m).
\]
If the inequality is strict, the first coordinate already gives
$\Lambda_m(q)<_{\mathrm{lex}}\Lambda_m(p_m)$.

Suppose henceforth that $q(a_m)=1/(m+1)$. Since its initial weight is
zero, at least one operation contains $a_m$. Define
\[
t_m=\max\{s\in\{1,\ldots,k\}:a_m\in S_s\},
\]
the time immediately after the last operation containing $a_m$.
By the equality conditions in the root-source barrier,
$S_{t_m}=P_m$ and
\[
q^{(t_m)}(v)=\frac1{m+1}\quad(v\in P_m),
\qquad
q^{(t_m)}(b)=q^{(t_m)}(c)=0.
\]
In particular, $t_m\ge1$. No later operation contains $a_m$, so every
later selected set is contained in $P_{m-1}\cup\{b,c\}$.

We next compare the coordinate at $a_{m-1}$. For $t_m\le s\le k$ and
$v\in P_{m-1}\cup\{b,c\}$, define
\[
z^{(s)}(v)=1-(m+1)q^{(s)}(v).
\]
At time $t_m$,
\[
z^{(t_m)}(v)=0\quad(v\in P_{m-1}),
\qquad
z^{(t_m)}(b)=z^{(t_m)}(c)=1.
\]
By \eqref{eq:affine-average}, the subsequent vectors $z^{(s)}$ follow
the same averaging sequence, and every operation remains effective.
Applying the barrier for two leaf sources on $P_{m-1}\cup\{b,c\}$ gives
\[
z^{(k)}(a_{m-1})\le\frac{2m+1}{(m+1)^2}.
\]
Transforming back to $q(a_{m-1})$, we obtain
\[
\begin{aligned}
q(a_{m-1})
&=\frac{1-z^{(k)}(a_{m-1})}{m+1}\\
&\ge\frac1{m+1}
   \left(1-\frac{2m+1}{(m+1)^2}\right)\\
&=\frac{m^2}{(m+1)^3}
=\rho_{m-1}=p_m(a_{m-1}).
\end{aligned}
\]
Hence $-q(a_{m-1})\le-p_m(a_{m-1})$. If this is strict, the first
coordinates agree and the second coordinate is strictly smaller, again
giving $\Lambda_m(q)<_{\mathrm{lex}}\Lambda_m(p_m)$.

Now suppose that $q(a_{m-1})=\rho_{m-1}$, so equality holds in the
barrier for two leaf sources. Since $z^{(t_m)}(a_{m-1})=0$, whereas its
final value is $(2m+1)/(m+1)^2>0$, some operation after $t_m$ contains
$a_{m-1}$. Define
\[
t_{m-1}=\max\{s\in\{1,\ldots,k\}:a_{m-1}\in S_s\}.
\]
Apply the barrier for two leaf sources up to time $t_{m-1}$.
Its equality conditions give $t_{m-1}-t_m\ge2$. At time $t_{m-1}$,
all transformed path weights and one transformed leaf weight equal
$(2m+1)/(m+1)^2$, while the other leaf has transformed weight $1/(m+1)$.
Transforming back gives
\[
q^{(t_{m-1})}(v)=\rho_{m-1}\quad(v\in P_{m-1}),
\qquad
\{q^{(t_{m-1})}(b),q^{(t_{m-1})}(c)\}
=\{\rho_{m-1},\rho_m\}.
\]
The weight at the other leaf is
$\bigl(1-1/(m+1)\bigr)/(m+1)=\rho_m$.
All later selected sets avoid $a_{m-1}$ and $a_m$.

We continue with $a_{m-2},\ldots,a_0$. Suppose that $1\le j\le m-1$
and the preceding comparisons have given
\[
q(a_m)=p_m(a_m),\qquad
q(a_i)=p_m(a_i)\quad(j\le i\le m-1),
\]
together with the time $t_j$ immediately after the last operation
containing $a_j$. Our induction hypothesis is that, at time $t_j$,
\begin{equation}\label{eq:equality-state}
q^{(t_j)}(v)=\rho_j\quad(v\in P_j),
\qquad
\{q^{(t_j)}(b),q^{(t_j)}(c)\}=\{\rho_j,\rho_{j+1}\},
\end{equation}
and all later selected sets avoid $a_j,a_{j+1},\ldots,a_m$.
The preceding calculation establishes this hypothesis for $j=m-1$.

For $t_j\le s\le k$ and $v\in P_{j-1}\cup\{b,c\}$, redefine the
transformed weights for the present comparison by
\[
z^{(s)}(v)
=\frac{q^{(s)}(v)-\rho_j}{\rho_{j+1}-\rho_j}.
\]
The denominator is nonzero by \eqref{eq:reference-signs}.
By \eqref{eq:equality-state}, at time $t_j$ every coordinate on $P_{j-1}$
is zero, and the coordinates at $b,c$ are zero and one in some order.
This transformation is illustrated in Figure~\ref{fig:affine-source}.

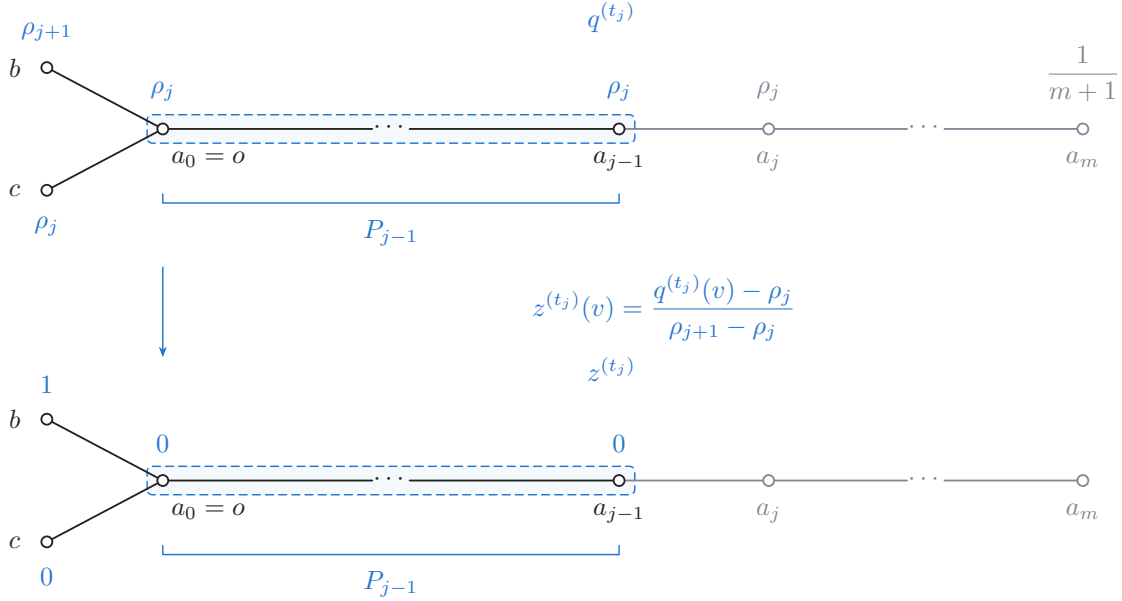
\begin{figure}[!htbp]
\centering
\resizebox{.96\linewidth}{!}{
\begin{tikzpicture}[wt/base]
\path (0,-6.5) rectangle (15.4,1.8);
\draw[wt/selection] (1.990,-0.190) rectangle (8.510,0.190);
\draw[wt/selection] (1.990,-4.890) rectangle (8.510,-4.510);
\draw[wt/edge] (0.650,0.820) -- (2.200,0.000);
\draw[wt/edge] (0.650,-0.820) -- (2.200,0.000);
\draw[wt/edge] (2.200,0.000) -- (5.000,0.000);
\draw[wt/edge] (5.500,0.000) -- (8.300,0.000);
\draw[wt/ghost] (8.300,0.000) -- (10.300,0.000);
\draw[wt/ghost] (10.300,0.000) -- (12.150,0.000);
\draw[wt/ghost] (12.650,0.000) -- (14.500,0.000);
\draw[wt/bracket,draw=wtblue] (2.200,-0.870) -- (2.200,-0.990) -- (8.300,-0.990) -- (8.300,-0.870);
\draw[wt/edge] (0.650,-3.880) -- (2.200,-4.700);
\draw[wt/edge] (0.650,-5.520) -- (2.200,-4.700);
\draw[wt/edge] (2.200,-4.700) -- (5.000,-4.700);
\draw[wt/edge] (5.500,-4.700) -- (8.300,-4.700);
\draw[wt/ghost] (8.300,-4.700) -- (10.300,-4.700);
\draw[wt/ghost] (10.300,-4.700) -- (12.150,-4.700);
\draw[wt/ghost] (12.650,-4.700) -- (14.500,-4.700);
\draw[wt/bracket,draw=wtblue] (2.200,-5.570) -- (2.200,-5.690) -- (8.300,-5.690) -- (8.300,-5.570);
\draw[wt/bracket,-{Stealth[length=4pt]}] (2.2,-1.85) -- (2.2,-3.05);
\draw[draw=wtblack,fill=white,line width=.70pt] (2.200,0.000) circle[radius=.073cm];
\draw[draw=wtblack,fill=white,line width=.70pt] (8.300,0.000) circle[radius=.073cm];
\draw[draw=wtgray,fill=white,line width=.70pt] (10.300,0.000) circle[radius=.073cm];
\draw[draw=wtgray,fill=white,line width=.70pt] (14.500,0.000) circle[radius=.073cm];
\draw[wt/edge,fill=white] (0.650,0.820) circle[radius=.073cm];
\draw[wt/edge,fill=white] (0.650,-0.820) circle[radius=.073cm];
\draw[draw=wtblack,fill=white,line width=.70pt] (2.200,-4.700) circle[radius=.073cm];
\draw[draw=wtblack,fill=white,line width=.70pt] (8.300,-4.700) circle[radius=.073cm];
\draw[draw=wtgray,fill=white,line width=.70pt] (10.300,-4.700) circle[radius=.073cm];
\draw[draw=wtgray,fill=white,line width=.70pt] (14.500,-4.700) circle[radius=.073cm];
\draw[wt/edge,fill=white] (0.650,-3.880) circle[radius=.073cm];
\draw[wt/edge,fill=white] (0.650,-5.520) circle[radius=.073cm];
\node[wt/math,text=wtblack,anchor=center] at (5.250,0.018) {$\displaystyle \cdots$};
\node[wt/math,text=wtgray,anchor=center] at (12.400,0.018) {$\displaystyle \cdots$};
\node[wt/math,text=wtblue,anchor=south] at (2.200,0.340) {$\displaystyle \rho_j$};
\node[wt/math,text=wtblack,anchor=north west] at (2.280,-0.270) {$\displaystyle a_0=o$};
\node[wt/math,text=wtblack,anchor=north] at (8.300,-0.270) {$\displaystyle a_{j-1}$};
\node[wt/math,text=wtblue,anchor=south] at (8.300,0.340) {$\displaystyle \rho_j$};
\node[wt/math,text=wtgray,anchor=north] at (10.300,-0.270) {$\displaystyle a_j$};
\node[wt/math,text=wtgray,anchor=south] at (10.300,0.340) {$\displaystyle \rho_j$};
\node[wt/math,text=wtgray,anchor=north] at (14.500,-0.270) {$\displaystyle a_m$};
\node[wt/math,text=wtgray,anchor=south] at (14.500,0.340) {$\displaystyle \frac1{m+1}$};
\node[wt/math,text=wtblack,anchor=east] at (0.330,0.820) {$\displaystyle b$};
\node[wt/math,text=wtblue,anchor=south] at (0.650,1.140) {$\displaystyle \rho_{j+1}$};
\node[wt/math,text=wtblack,anchor=east] at (0.330,-0.820) {$\displaystyle c$};
\node[wt/math,text=wtblue,anchor=north] at (0.650,-1.140) {$\displaystyle \rho_j$};
\node[wt/math,text=wtblue,anchor=center] at (8.200,1.490) {$\displaystyle q^{(t_j)}$};
\node[wt/math,text=wtblue,anchor=north] at (5.250,-1.180) {$\displaystyle P_{j-1}$};
\node[wt/math,text=wtblack,anchor=center] at (5.250,-4.682) {$\displaystyle \cdots$};
\node[wt/math,text=wtgray,anchor=center] at (12.400,-4.682) {$\displaystyle \cdots$};
\node[wt/math,text=wtblue,anchor=south] at (2.200,-4.360) {$\displaystyle 0$};
\node[wt/math,text=wtblack,anchor=north west] at (2.280,-4.970) {$\displaystyle a_0=o$};
\node[wt/math,text=wtblack,anchor=north] at (8.300,-4.970) {$\displaystyle a_{j-1}$};
\node[wt/math,text=wtblue,anchor=south] at (8.300,-4.360) {$\displaystyle 0$};
\node[wt/math,text=wtgray,anchor=north] at (10.300,-4.970) {$\displaystyle a_j$};
\node[wt/math,text=wtgray,anchor=north] at (14.500,-4.970) {$\displaystyle a_m$};
\node[wt/math,text=wtblack,anchor=east] at (0.330,-3.880) {$\displaystyle b$};
\node[wt/math,text=wtblue,anchor=south] at (0.650,-3.560) {$\displaystyle 1$};
\node[wt/math,text=wtblack,anchor=east] at (0.330,-5.520) {$\displaystyle c$};
\node[wt/math,text=wtblue,anchor=north] at (0.650,-5.840) {$\displaystyle 0$};
\node[wt/math,text=wtblue,anchor=center] at (8.200,-3.210) {$\displaystyle z^{(t_j)}$};
\node[wt/math,text=wtblue,anchor=north] at (5.250,-5.880) {$\displaystyle P_{j-1}$};
\node[wt/math,text=wtblue,anchor=center] at (8.900,-2.390) {$\displaystyle z^{(t_j)}(v)=\frac{q^{(t_j)}(v)-\rho_j}{\rho_{j+1}-\rho_j}$};
\end{tikzpicture}
}
\caption{From an equality configuration to a single leaf source, for $j\ge2$
and $q^{(t_j)}(b)=\rho_{j+1}$.}
\label{fig:affine-source}
\end{figure}

By \eqref{eq:affine-average}, the transformed weights continue to follow
the original averaging sequence. We may therefore apply the barrier for
a single leaf source, interchanging $b$ and $c$ if the initial unit
weight is at $c$. This barrier, including the boundary cases covered
by the parameter ranges specified in Section~2, gives
\[
z^{(k)}(a_{j-1})\le\frac1{j+1}.
\]
On the other hand, \eqref{eq:reference-recursion} can be written as
\[
\rho_{j-1}
=\rho_j+\frac{\rho_{j+1}-\rho_j}{j+1}.
\]
Combining this with
$q(a_{j-1})=\rho_j+(\rho_{j+1}-\rho_j)z^{(k)}(a_{j-1})$ gives
\[
\begin{cases}
q(a_{j-1})\le\rho_{j-1}, & \rho_{j+1}-\rho_j>0,\\
q(a_{j-1})\ge\rho_{j-1}, & \rho_{j+1}-\rho_j<0.
\end{cases}
\]
By \eqref{eq:reference-signs} and
$(-1)^{m-1-j}=(-1)^{m-j+1}$,
$(-1)^{m-j+1}(\rho_{j+1}-\rho_j)>0$.
The sign of the coordinate corresponding to $a_{j-1}$ in $\Lambda_m$
is precisely $(-1)^{m-j+1}$. Thus both cases give
\[
(-1)^{m-j+1}q(a_{j-1})
\le(-1)^{m-j+1}p_m(a_{j-1}).
\]
If this inequality is strict, all preceding coordinates agree and the
present coordinate is strictly smaller. Hence
$\Lambda_m(q)<_{\mathrm{lex}}\Lambda_m(p_m)$.

If the present coordinate also agrees, then
$z^{(k)}(a_{j-1})=1/(j+1)$. Its value at time $t_j$ is zero, so some
effective operation after $t_j$ contains $a_{j-1}$. Define
\[
t_{j-1}=\max\{s\in\{1,\ldots,k\}:a_{j-1}\in S_s\}.
\]
Then $t_{j-1}>t_j$. By the equality conditions in the barrier for a
single leaf source, this operation averages $P_{j-1}$ together with the
leaf whose transformed weight at time $t_j$ is one.
At time $t_{j-1}$, every transformed coordinate on this set equals
$1/(j+1)$, and the transformed weight at the other leaf is zero.
Transforming back gives
\[
q^{(t_{j-1})}(v)=\rho_{j-1}\quad(v\in P_{j-1}),
\qquad
\{q^{(t_{j-1})}(b),q^{(t_{j-1})}(c)\}
=\{\rho_{j-1},\rho_j\}.
\]
All later selected sets avoid $a_{j-1},a_j,\ldots,a_m$.
This is the configuration in \eqref{eq:equality-state} with $j$ replaced
by $j-1$, so the induction continues.

The comparisons cover $a_m,a_{m-1},\ldots,a_0$. Therefore,
$\Lambda_m(q)\le_{\mathrm{lex}}\Lambda_m(p_m)$ for every
$q\in\mathcal R_{\mathrm{fin}}$.
If all compared coordinates agree, then all times $t_m,\ldots,t_0$
are defined and satisfy
\[
t_m\ge1,\qquad
t_{m-1}-t_m\ge2,\qquad
t_{j-1}-t_j\ge1\quad(1\le j\le m-1).
\]
Consequently,
\[
\begin{aligned}
\ell(q)=k
&\ge t_0\\
&=t_m+(t_{m-1}-t_m)
  +\sum_{j=1}^{m-1}(t_{j-1}-t_j)\\
&\ge1+2+(m-1)=m+2.
\end{aligned}
\]
The operation counts forced by these equality times are shown in
Figure~\ref{fig:equality-times}.

\begin{figure}[!htbp]
\centering
\resizebox{.96\linewidth}{!}{
\begin{tikzpicture}[wt/base]
\path (0,-2.7) rectangle (15.4,2.2);
\draw[wt/edge] (0.800,0.000) -- (3.200,0.000);
\draw[wt/edge] (3.200,0.000) -- (6.200,0.000);
\draw[wt/edge] (6.200,0.000) -- (10.100,0.000);
\draw[wt/edge] (10.600,0.000) -- (14.500,0.000);
\draw[wt/bracket,draw=wtblue] (0.800,0.860) -- (0.800,0.980) -- (3.200,0.980) -- (3.200,0.860);
\draw[wt/bracket,draw=wtblue] (3.200,0.860) -- (3.200,0.980) -- (6.200,0.980) -- (6.200,0.860);
\draw[wt/bracket,draw=wtblue] (6.200,0.860) -- (6.200,0.980) -- (14.500,0.980) -- (14.500,0.860);
\draw[draw=wtblack,fill=white,line width=.70pt] (0.800,0.000) circle[radius=.073cm];
\draw[draw=wtblack,fill=white,line width=.70pt] (3.200,0.000) circle[radius=.073cm];
\draw[draw=wtblack,fill=white,line width=.70pt] (6.200,0.000) circle[radius=.073cm];
\draw[draw=wtblack,fill=white,line width=.70pt] (14.500,0.000) circle[radius=.073cm];
\node[wt/math,text=wtblue,anchor=center] at (7.700,1.960) {$\displaystyle \Lambda_m(q)=\Lambda_m(p_m)$};
\node[wt/math,text=wtblack,anchor=center] at (10.350,0.018) {$\displaystyle \cdots$};
\node[wt/math,text=wtblack,anchor=north] at (0.800,-0.270) {$\displaystyle 0$};
\node[wt/math,text=wtblack,anchor=north] at (3.200,-0.270) {$\displaystyle t_m$};
\node[wt/math,text=wtblack,anchor=north] at (3.200,-0.950) {$\displaystyle a_m$};
\node[wt/math,text=wtblack,anchor=north] at (6.200,-0.270) {$\displaystyle t_{m-1}$};
\node[wt/math,text=wtblack,anchor=north] at (6.200,-0.950) {$\displaystyle a_{m-1}$};
\node[wt/math,text=wtblack,anchor=north] at (14.500,-0.270) {$\displaystyle t_0$};
\node[wt/math,text=wtblack,anchor=north] at (14.500,-0.950) {$\displaystyle a_0$};
\node[wt/math,text=wtblue,anchor=south] at (2.000,1.170) {$\displaystyle \ge1$};
\node[wt/math,text=wtblue,anchor=south] at (4.700,1.170) {$\displaystyle \ge2$};
\node[wt/math,text=wtblue,anchor=south] at (10.350,1.170) {$\displaystyle \ge m-1$};
\node[wt/math,text=wtblue,anchor=center,font=\footnotesize] at (10.100,-1.210) {$\displaystyle t_{j-1}-t_j\ge1\quad(1\le j\le m-1)$};
\node[wt/math,text=wtblue,anchor=center] at (7.700,-2.240) {$\displaystyle \ell(q)\ge t_0\ge1+2+(m-1)=m+2$};
\end{tikzpicture}
}
\caption{Effective operation counts when all compared coordinates agree.}
\label{fig:equality-times}
\end{figure}
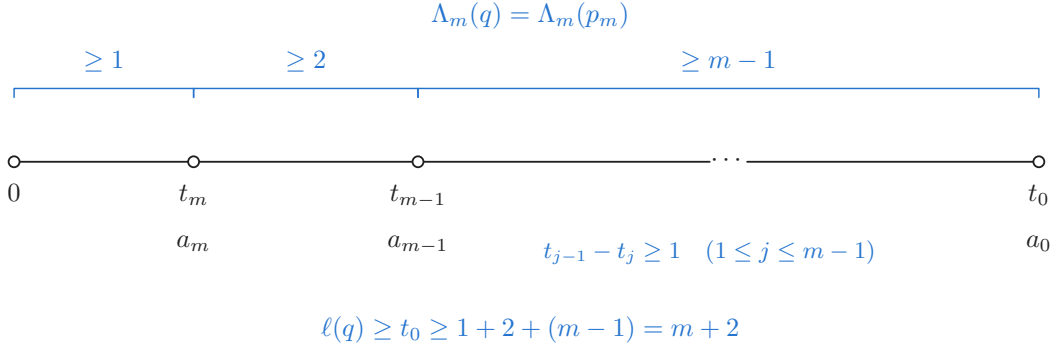

Finally, if $q\in\mathcal R_{m+1}$, then $\ell(q)\le m+1$.
Thus not all compared coordinates can agree, and necessarily
$\Lambda_m(q)<_{\mathrm{lex}}\Lambda_m(p_m)$.
\end{proof}

In particular, taking $q=p_m$ gives $\ell(p_m)\ge m+2$.
Together with the construction using $m+2$ operations in Section~3,
this yields $\ell(p_m)=m+2$. In the next section, we choose nonnegative
initial weights $\omega_m$ that turn the strict lexicographic comparison
in Theorem~\ref{thm:profile-comparison} into a strict inner product inequality.

\FloatBarrier
\section{Initial weights and inner product comparison}

We construct $\omega_m$ so that $\langle p_m-q,\omega_m\rangle>0$
for every $q\in\mathcal R_{m+1}$. We first bound the nonzero coordinate
differences between $p_m$ and these profiles.

Starting from $e_o$, suppose that immediately before an operation all
coordinates can be written as fractions with a common positive integer
denominator $D$. After averaging a set $S$, they can all be written
with denominator $D|S|$. Hence, after $r$ operations, the reduced
denominator of each coordinate is at most the product of the sizes of
the $r$ selected sets, and therefore at most $(m+3)^r$.
Here we use $|V(G)|=m+3$; when $r=0$, the empty product is one.

The profile $p_m$ is obtained in $m+2$ operations, whereas
$q\in\mathcal R_{m+1}$ requires at most $m+1$ operations.
Their coordinate denominators in lowest terms are therefore at most
$(m+3)^{m+2}$ and $(m+3)^{m+1}$, respectively.
Put $\delta=(m+3)^{-2(m+2)}$.
The absolute difference between distinct fractions is at least the
reciprocal of the product of their denominators. Thus, for every
$q\in\mathcal R_{m+1}$ and $v\in V(G)$,
\begin{equation}\label{eq:coordinate-gap}
p_m(v)\ne q(v)
\quad\Longrightarrow\quad
|p_m(v)-q(v)|
\ge(m+3)^{-(2m+3)}>\delta.
\end{equation}

In the initial weights below, the signed term at $a_i$ will be
$(-1)^{m-i}Q^i$. We choose $Q$ so that $1/(Q-1)<\delta$.
This ensures that the first unequal coordinate in lexicographic order
determines the sign of the inner product difference.
We then add a constant $B$ to the initial weight at every vertex.
The largest absolute value among the negative terms is $Q^{m-1}$,
so $B=Q^{m-1}$ makes all initial weights nonnegative.
Accordingly, choose $Q=4(m+3)^{2(m+2)}$ and $B=Q^{m-1}$.
Both are positive integers, and
\[
\frac1{Q-1}=\frac{\delta}{4-\delta}<\delta,
\]
where we use $0<\delta<1$.

\begin{theorem}[Strict comparison for nonnegative integer weights]
\label{thm:weighted-comparison}
Define initial weights $\omega_m$ on $G=S(m,1,1)$ by
\begin{equation}\label{eq:initial-weights}
\omega_m(a_i)=B+(-1)^{m-i}Q^i\quad(0\le i\le m),
\qquad
\omega_m(b)=\omega_m(c)=B.
\end{equation}
Then $\omega_m$ is a nonnegative integer-valued weight function and,
for every $q\in\mathcal R_{m+1}$,
\[
\langle q,\omega_m\rangle<\langle p_m,\omega_m\rangle.
\]
\end{theorem}

\begin{proof}
All weights in \eqref{eq:initial-weights} are integers.
If $(-1)^{m-i}=1$, then $\omega_m(a_i)>0$.
If $(-1)^{m-i}=-1$, then $i\le m-1$, so
$Q^i\le Q^{m-1}=B$ and $\omega_m(a_i)=B-Q^i\ge0$.
Also, $\omega_m(b)=\omega_m(c)=B>0$.
Thus $\omega_m$ is a nonnegative integer-valued weight function.

Fix $q\in\mathcal R_{m+1}$. By
Theorem~\ref{thm:profile-comparison},
$\Lambda_m(q)<_{\mathrm{lex}}\Lambda_m(p_m)$.
Define
\[
s=\max\{i\in\{0,\ldots,m\}:p_m(a_i)\ne q(a_i)\}.
\]
Thus $a_s$ is the first unequal coordinate in the order
$a_m,a_{m-1},\ldots,a_0$. The strict lexicographic inequality ensures
that the defining set is nonempty and gives
\[
p_m(a_i)=q(a_i)\quad(s<i\le m),
\qquad
(-1)^{m-s}\bigl(p_m(a_s)-q(a_s)\bigr)>0.
\]
By \eqref{eq:coordinate-gap}, we further have
\[
(-1)^{m-s}\bigl(p_m(a_s)-q(a_s)\bigr)>\delta.
\]
Since $p_m(v),q(v)\in[0,1]$, we also have
$|p_m(v)-q(v)|\le1$ for every vertex.

Since the coordinates of both profiles sum to one, the constant part
of \eqref{eq:initial-weights} contributes
\[
B\sum_{v\in V(G)}\bigl(p_m(v)-q(v)\bigr)=B(1-1)=0.
\]
Therefore,
\[
\begin{aligned}
\langle p_m-q,\omega_m\rangle
&=\sum_{i=0}^{s}(-1)^{m-i}Q^i
   \bigl(p_m(a_i)-q(a_i)\bigr)\\
&>Q^s\delta-\sum_{i=0}^{s-1}Q^i\\
&=Q^s\delta-\frac{Q^s-1}{Q-1}\\
&>Q^s\left(\delta-\frac1{Q-1}\right)>0.
\end{aligned}
\]
When $s=0$, the empty sum is zero and the same calculation applies.
Thus $\langle q,\omega_m\rangle<\langle p_m,\omega_m\rangle$.
\end{proof}

\FloatBarrier
\section{A lower bound on the length of a shortest optimal strategy}

\begin{theorem}[Lower bound on the length of a shortest optimal strategy]
\label{thm:main}
For every integer $m\ge2$, let $\omega_m$ be the nonnegative integer
initial weight function defined in \eqref{eq:initial-weights}. Then
\[
t\bigl(S(m,1,1),o,\omega_m\bigr)\ge m+2.
\]
\end{theorem}

\begin{proof}
Apply the $m+2$ operations constructing $p_m$ in Section~3 to the
initial weights $\omega_m$ in reverse order. The reverse profile of
this forward strategy is $p_m$, so \eqref{eq:forward-reverse} shows
that its resulting target weight at $o$ is
$\langle p_m,\omega_m\rangle$.

On the other hand, the reverse profile $q$ of any forward strategy
using at most $m+1$ operations belongs to $\mathcal R_{m+1}$.
By Theorem~\ref{thm:weighted-comparison}, the target weight produced
by such a strategy satisfies
\[
\langle q,\omega_m\rangle<\langle p_m,\omega_m\rangle.
\]
The right-hand side is itself reachable. Hence no strategy using at
most $m+1$ operations can attain the maximum reachable target weight.
Every optimal strategy therefore requires at least $m+2$ operations, that is,
$t\bigl(S(m,1,1),o,\omega_m\bigr)\ge m+2$.
\end{proof}

For any positive integer $K$, choose an integer
$m\ge\max\{2,K\}$. The theorem gives
\[
t\bigl(S(m,1,1),o,\omega_m\bigr)\ge m+2>K.
\]
Thus the lengths of shortest optimal strategies on $3$-spiders have
no uniform finite bound. Since every graph in this family has maximum
degree three, the same conclusion holds even when the maximum degree
is fixed at three.

\section*{Acknowledgements}
The authors are grateful to Hehui Wu and Ningyuan Yang for helpful discussions and valuable suggestions during the early stages of this project. The authors also thank ChatGPT (GPT-5.6 Sol) for assistance with translation, typesetting, writing, and merging steps. The mathematical proofs in this paper were almost entirely carried out by the human authors.


\clearpage
\begin{thebibliography}{9}
\bibitem{gollin2024}
J.~P.~Gollin, K.~Hendrey, H.~Huang, T.~Huynh, B.~Mohar, S.-i.~Oum,
N.~Yang, W.-H.~Yu and X.~Zhu,
\emph{Sharing tea on a graph},
\href{https://arxiv.org/abs/2405.15353}{arXiv:2405.15353}, 2024.

\bibitem{tao2025}
T.~Tao,
\emph{Water transport on a path: Finding the strategy through its existence},
Advances and Applications in Discrete Mathematics \textbf{42} (2025),
no.~5, 471--480.
\href{https://doi.org/10.17654/0974165825031}{doi:10.17654/0974165825031}.

\bibitem{haggstrom2012}
O.~H\"aggstr\"om,
\emph{A pairwise averaging procedure with application to consensus formation
in the Deffuant model},
Acta Applicandae Mathematicae \textbf{119} (2012), no.~1, 185--201.
\href{https://doi.org/10.1007/s10440-011-9668-9}{doi:10.1007/s10440-011-9668-9}.

\bibitem{haggstrom2019}
O.~H\"aggstr\"om and T.~Hirscher,
\emph{Water transport on infinite graphs},
Random Structures \& Algorithms \textbf{54} (2019), no.~3, 515--527.
\href{https://doi.org/10.1002/rsa.20801}{doi:10.1002/rsa.20801}.

\bibitem{vilkas2025}
T.~Vilkas,
\emph{Water transport on finite graphs},
\href{https://arxiv.org/abs/2501.16911}{arXiv:2501.16911}, 2025.
\end{thebibliography}
\end{document}